%% file: main_v1.tex
\documentclass[11pt,a4paper]{article}
\input{settings/packages.tex}
\input{settings/style}

\input{settings/macros}

\title{\textsc{Convergence Rates for the MAP Estimator in PDE-Regression Models with Random Design}}
\author{\textsc{Maximilian Siebel}\thanks{\texttt{siebel@math.uni-heidelberg.de}}}
\affil{\textsc{Institute for Mathematics, Heidelberg University}}
\date{}

\begin{document}

\maketitle

\begin{abstract}
    We consider the statistical inverse problem of recovering a parameter $\theta\in H^\alpha$ from data arising from the Gaussian regression problem 
    \begin{equation*}
        Y = \scrG(\theta)(Z)+\varepsilon
    \end{equation*}
    with nonlinear forward map $\scrG:\IL^2\to\IL^2$, random design point $Z$ and Gaussian noise $\varepsilon$. The estimation strategy is based on a least squares approach under $\norm{\cdot}_{H^\alpha}$-constraints. We establish the existence of a least squares estimator $\hat{\theta}$ as a maximizer for a given functional under Lipschitz-type assumptions on the forward map $\scrG$.  A general concentration result is shown, which is used to prove consistency and upper bounds for the prediction error. The corresponding rates of convergence reflect not only the smoothness of the parameter of interest but also the ill-posedness of the underlying inverse problem. We apply the general model to the Darcy problem, where the recovery of an unknown coefficient function $f$ of a PDE is of interest. For this example, we also provide corresponding rates of convergence for the prediction and estimation errors. Additionally, we briefly discuss the applicability of the general model to other problems.
\end{abstract}

\tableofcontents

\input{content/intro}

\input{content/model/model}

\input{content/pde/pde}

\newpage
\section{Appendix}\label{sec:app}
\input{content/appendix/app_map.tex}

\input{content/appendix/app_pde.tex}
\input{content/appendix/auxres.tex}

\section*{Acknowledgements}
The author thanks Richard Nickl and Jan Johannes for their time, patience, and fruitful discussions, as well as for the visit to
Richard Nickl’s research group in Cambridge in November 2023, which marked the beginning of this work.
This work is funded by Deutsche Forschungsgemeinschaft (DFG, German Re-
search Foundation) under Germany’s Excellence Strategy EXC-2181/1-39090098 (the Heidelberg
STRUCTURES Cluster of Excellence).
\newpage
\appendix

\bibliographystyle{plainnat}
\bibliography{mapbib.bib}

\end{document}

%% file: settings/packages.tex
\usepackage{graphicx}
\usepackage[margin=1in]{geometry}
\usepackage{amsthm,amssymb}
\usepackage{bbm}
\usepackage[dvipsnames]{xcolor}
\usepackage{amsfonts}
\usepackage{amsmath}
\usepackage{mathtools}
\usepackage{dsfont}
\usepackage{enumerate}   
\usepackage{stmaryrd}
\usepackage[toc,page]{appendix}
\usepackage[colorlinks=true, citecolor=blue, linkcolor=black, urlcolor=blue]{hyperref}
\usepackage{graphicx,float}
\usepackage{natbib}  
\usepackage{mdframed}
\usepackage{soul}
\usepackage{mathrsfs} 
\usepackage{cleveref}
\usepackage{authblk}
\usepackage{amsthm}
\usepackage{thmtools}
\usepackage{enumitem}

%% file: settings/style.tex
\definecolor{cardinalred}{RGB}{140,21,21}
\crefformat{definition}{\textcolor{Blue}{Definition}~#2{\color{Blue}\textcolor{gray}{\S\phantom{.}#1}}#3}
\crefformat{theorem}{\textcolor{cardinalred}{Theorem}~#2{\color{cardinalred}\textcolor{gray}{\S\phantom{.}#1}}#3}
\crefformat{assumption}{\textcolor{cardinalred}{Condition}~#2{\color{cardinalred}\textcolor{gray}{\S\phantom{.}#1}}#3}
\crefformat{corollary}{\textcolor{cardinalred}{Corollary}~#2{\color{cardinalred}\textcolor{gray}{\S\phantom{.}#1}}#3}
\crefformat{proposition}{\textcolor{cardinalred}{Proposition}~#2{\color{cardinalred}\textcolor{gray}{\S\phantom{.}#1}}#3}
\crefformat{lemma}{\textcolor{cardinalred}{Lemma}~#2{\color{cardinalred}\textcolor{gray}{\S\phantom{.}#1}}#3}
\crefformat{example}{\textcolor{OliveGreen}{Example}~#2{\color{OliveGreen}\textcolor{gray}{\S\phantom{.}#1}}#3}
\crefformat{notation}{\textcolor{oldgold}{Notation}~#2{\color{oldgold}\textcolor{gray}{\S\phantom{.}#1}}#3}
\crefformat{remark}{\textcolor{Brown}{Remark}~#2{\color{Brown}\textcolor{gray}{\S\phantom{.}#1}}#3}
\crefformat{appendix}{\textcolor{Black}{\textsc{Appendix}}~#2{\color{Black}\textcolor{gray}{\S\phantom{.}#1}}#3}

\crefformat{section}{\textcolor{Black}{\textsc{Section}}~#2{\color{Black}\textcolor{gray}{\S\phantom{.}#1}}#3}
\crefformat{equation}{\textcolor{Black}{\textsc{Equation}}~#2{\color{Black}\textcolor{gray}{#1}}#3}

\newcommand{\magenta}[1]{\textcolor{black}{#1}}
\definecolor{bananayellow}{rgb}{1.0, 0.88, 0.21}
\definecolor{citrine}{rgb}{0.89, 0.82, 0.04}
\definecolor{goldenyellow}{rgb}{1.0, 0.87, 0.0}
\definecolor{jonquil}{rgb}{0.98, 0.85, 0.37}
\definecolor{mikadoyellow}{rgb}{1.0, 0.77, 0.05}
\definecolor{mustard}{rgb}{1.0, 0.86, 0.35}
\definecolor{oldgold}{rgb}{0.81, 0.71, 0.23}

\declaretheoremstyle[
qed=\raisebox{1pt}{\textcolor{cardinalred}{\openbox}},
headfont = \color{cardinalred},
headformat = \NAME{} \textcolor{gray}{\NUMBER}:\NOTE,
]{theoremstyles}
\declaretheorem[
style=theoremstyles,
name=\textsc{Theorem \textcolor{gray}{\S}},
numberwithin=section,
]{theorem}

\declaretheoremstyle[
qed=\raisebox{1pt}{\textcolor{cardinalred}{\openbox}},
headfont = \color{cardinalred},
headformat = \NAME{} \textcolor{gray}{\NUMBER}:\NOTE,
]{Corollarystyles}
\declaretheorem[
style=Corollarystyles,
name=\textsc{Corollary \textcolor{gray}{\S}},
numberlike=theorem,
]{corollary}

\declaretheoremstyle[
qed=\raisebox{1pt}{\textcolor{cardinalred}{\openbox}},
headfont = \color{cardinalred},
headformat = \NAME{} \textcolor{gray}{\NUMBER}:\NOTE,
]{Propostyles}
\declaretheorem[
style=Propostyles,
name=\textsc{Proposition \textcolor{gray}{\S}},
numberlike=theorem,
]{proposition}

\declaretheoremstyle[
qed=\raisebox{1pt}{\textcolor{cardinalred}{\openbox}},
headfont = \color{cardinalred},
headformat = \NAME{} \textcolor{gray}{\NUMBER}:\NOTE,
]{Lemmastyles}
\declaretheorem[
style=Lemmastyles,
name=\textsc{Lemma \textcolor{gray}{\S}},
numberlike=theorem,
]{lemma}

\declaretheoremstyle[
qed=\raisebox{1pt}{\textcolor{Blue}{\openbox}},
headfont = \protect\color{Blue},
headformat = \NAME{} \textcolor{gray}{\NUMBER}:\NOTE,
]{definitions}
\declaretheorem[
style=definitions,
name=\textsc{Definition \textcolor{gray}{\S}},
numberlike=theorem,
]{definition}

\declaretheoremstyle[
qed=\raisebox{1pt}{\textcolor{OliveGreen}{\openbox}},
headfont = \color{OliveGreen},
headformat = \NAME{} \textcolor{gray}{\NUMBER}:\NOTE,
]{examples}

\declaretheoremstyle[
qed=\raisebox{1pt}{\textcolor{oldgold}{\openbox}},
headfont = \color{oldgold},
headformat = \NAME{} \textcolor{gray}{\NUMBER}:\NOTE,
]{notations}
\declaretheorem[
style=notations,
name=\textsc{Notation \textcolor{gray}{\S}},
numberlike=theorem,
]{notation}

\declaretheoremstyle[
qed=\raisebox{1pt}{\textcolor{Brown}{\openbox}},
headfont = \color{Brown},
headformat = \NAME{} \textcolor{gray}{\NUMBER}:\NOTE,
]{remarks}
\declaretheorem[
style=remarks,
name=\textsc{Remark \textcolor{gray}{\S}},
numberlike=theorem,
]{remark}

\declaretheoremstyle[
qed=\raisebox{1pt}{\textcolor{black}{\openbox}},
headfont = \color{cardinalred},
headformat = \NAME{} \textcolor{gray}{\NUMBER}:\NOTE,
]{assumptions}
\declaretheorem[
style=assumptions,
name=\textsc{Condition \textcolor{gray}{\S}},
numberlike=theorem,
]{assumption}

\declaretheoremstyle[
numbered=no,
notebraces={}{},
qed=\raisebox{1pt}{\textcolor{gray}{\openbox}},
headfont = \itshape\color{gray},
]{proofs}

%% file: settings/macros.tex
\newcommand{\IB}{\mathbb{B}}

\newcommand{\IE}{\mathbb{E}}

\newcommand{\IK}{\mathbb{K}}
\newcommand{\IL}{\mathbb{L}}

\newcommand{\IN}{\mathbb{N}}

\newcommand{\IP}{\mathbb{P}}

\newcommand{\IR}{\mathbb{R}}

\newcommand{\IZ}{\mathbb{Z}}

\newcommand{\calA}{\mathcal{A}}

\newcommand{\calC}{\mathcal{C}}

\newcommand{\calF}{\mathcal{F}}

\newcommand{\calH}{\mathcal{H}}

\newcommand{\calJ}{\mathcal{J}}

\newcommand{\calL}{\mathcal{L}}

\newcommand{\calO}{\mathcal{O}}

\newcommand{\calT}{\mathcal{T}}

\newcommand{\calW}{\mathcal{W}}

\newcommand{\calZ}{\mathcal{Z}}

\newcommand{\scrB}{\mathscr{B}}

\newcommand{\scrG}{\mathscr{G}}
\newcommand{\scrH}{\mathscr{H}}

\newcommand{\scrJ}{\mathscr{J}}

\newcommand{\scrZ}{\mathscr{Z}}

\newcommand{\frakd}{\mathfrak{d}}

\newcommand{\frakm}{\mathfrak{m}}

\newcommand{\frako}{\mathfrak{o}}

\newcommand{\frakD}{\mathfrak{D}}

\newcommand{\empM}{\widehat{\IP}_N}

\newcommand{\domain}{\calO}

\newcommand{\topo}{\calT}

\newcommand{\norm}[1]{\|#1\|}
\newcommand{\Vnorm}[1]{|#1|_V}

\newcommand{\Wnorm}[1]{\left|#1\right|_W}

\newcommand{\thetaMAP}{\hat{\theta}_{\textsc{MAP}}}
\newcommand{\fMAP}{\hat{f}_{\textsc{MAP}}}
\newcommand{\thetatrue}{\theta_\frako}
\newcommand{\thetatilde}{\tilde{\theta}}
\newcommand{\Thetatilde}{\tilde{\Theta}}

\newcommand{\nset}[1]{\left\llbracket #1\right\rrbracket}

\newcommand{\nnset}[2]{\left\llbracket #1,#2\right\rrbracket}

\newcommand{\const}{c}
\newcommand{\Const}{C}

\DeclareMathOperator*{\argmax}{arg\,max}

\newcommand{\Lip}{\operatorname{Lip}}

\newcommand{\dist}{\frakd}
\newcommand{\Dist}{\frakD}
\newcommand{\fmin}{f_{\min}}
\newcommand{\Ptt}{\texttt{P}}
\newcommand*{\ob}[1]{%
	\makebox[0pt][l]{%
		\kern.3em
		\raisebox{1.8ex}{\relsize{-3}{$\circ$}}
	}%
	#1%
}

\makeatletter
\newsavebox\myboxA
\newsavebox\myboxB
\newlength\mylenA

\newcommand*\xoverline[2][0.75]{%
	\sbox{\myboxA}{$\m@th#2$}%
	\setbox\myboxB\null% Phantom box
	\ht\myboxB=\ht\myboxA%
	\dp\myboxB=\dp\myboxA%
	\wd\myboxB=#1\wd\myboxA% Scale phantom
	\sbox\myboxB{$\m@th\overline{\copy\myboxB}$}%  Overlined phantom
	\setlength\mylenA{\the\wd\myboxA}%   calc width diff
	\addtolength\mylenA{-\the\wd\myboxB}%
	\ifdim\wd\myboxB<\wd\myboxA%
	\rlap{\hskip 0.5\mylenA\usebox\myboxB}{\usebox\myboxA}%
	\else
	\hskip -0.5\mylenA\rlap{\usebox\myboxA}{\hskip 0.5\mylenA\usebox\myboxB}%
	\fi}
\makeatother

\newcommand{\pRZ}{\mathbb{R}_{>0}}
\newcommand{\pRZz}{\mathbb{R}_{\geq0}}
\newcommand{\rate}{\mathfrak{r}}
\newcommand{\divergence}{\operatorname{div}}

\newcommand{\Zygmund}[2]{C_{\scalebox{0.6}{$\mathrm{Zy}$}}^{\mathrm{#1}}(#2)}

\newcommand{\floor}[1]{\lfloor #1 \rfloor}

\newcommand{\boldbeta}{\boldsymbol{\beta}}

\usepackage{etoolbox}
\makeatletter
\providecommand*{\cupdot}{%
	\mathbin{%
		\mathpalette\@cupdot{}%
	}%
}
\newcommand*{\@cupdot}[2]{%
	\ooalign{%
		$\m@th#1\cup$\cr
		\hidewidth$\m@th#1\cdot$\hidewidth
	}%
}
\makeatother
\makeatletter
\providerobustcmd*{\bigcupdot}{%
	\mathop{%
		\mathpalette\bigop@dot\bigcup
	}%
}
\newrobustcmd*{\bigop@dot}[2]{%
	\setbox0=\hbox{$\m@th#1#2$}%
	\vbox{%
		\lineskiplimit=\maxdimen
		\lineskip=-0.7\dimexpr\ht0+\dp0\relax
		\ialign{%
			\hfil##\hfil\cr
			$\m@th\cdot$\cr
			\box0\cr
		}%
	}%
}
\makeatother

\newcommand{\Hilbert}[1]{$(\mathbb{H},\langle\,\cdot\,|\,\cdot\,\rangle_{\mathbb{H}})$}
\newcommand{\Lz}[1]{\text{$\mathscr{L}^{2}([0,1])$}}

\def\DHLhksqrt#1#2{\setbox0=\hbox{$#1\sqrt{#2\,}$}\dimen0=\ht0
	\advance\dimen0-0.2\ht0
	\setbox2=\hbox{\vrule height\ht0 depth -\dimen0}%
	{\box0\lower0.4pt\box2}}

%% file: content/intro.tex
\section{Introduction}\label{sec:Intro}
 In this paper, we consider a nonlinear Gaussian regression model with random design, i.e., for a given sample size $N\in\IN$ and a fixed variance $\sigma^2\in\pRZ$, we have access to data $(Y_i,X_i)_{i=1}^N$ arising from
\begin{equation}\label{eq:IntroData}
	Y_i = u_f(X_i) + \sigma\varepsilon_i,\, i = 1,\dots, N,
\end{equation}
where $\varepsilon_1,\dots,\varepsilon_N$ are i.i.d. copies of $\varepsilon\sim N(0,1)$. Here, $X_1,\dots,X_N$ are i.i.d. copies of $X$, which is drawn independently of $\varepsilon$ and follows the uniform distribution over a bounded domain $\domain \subseteq \IR^d$ with smooth boundary $\partial \domain$. The regression function $u_f:\domain\to\IR$ denotes the solution of an (elliptic) partial differential equation (PDE), which depends nonlinearly on an \textit{unknown} coefficient function $f$. The goal is to present an estimation strategy for $f$ by using the data set $(Y_1,X_1),\dots,(Y_N,X_N)$, which is an example of a nonlinear statistical inverse problem as $f$ is not observed directly, see, e.g. \cite{Nickl_2023}.\\

Many real-world phenomena, like physical and biological laws, fit into this framework as they are often described by (non-)linear dynamical systems and partial differential equations which exhibit spatial and/or temporal dependencies. Examples range from the description of cell movements (see \cite{Tetsuji2021}) and geophysical laws describing the atmosphere (see \cite{Kalnay_2002}), to ocean dynamics (see \cite{Bennett2002}), turbulence (see \cite{Majda_Harlim_2012}), and fluid flows (see \cite{Cotter_2009}). Moreover, this setting aligns with various frameworks from the perspective of inverse problems: examples include electrical impedance tomography (see \cite{Calderon_2006} and \cite{Isaacson_2004}), photoacoustic tomography and other hybrid imaging problems (see \cite{Kuchment_2015,Bal_2014}, and particularly \cite{Bal_2010} for a general framework), inverse scattering (see \cite{Colton_2019} and \cite{Hohage_2015}), interacting particle models (see \cite{NicklGrigoriosRay_2024}), and reaction diffusion equations (see \cite{Nickl_2024}). In this work, we consider the so-called \textit{Darcy Problem} (see, for instance \cite{Stuart_2010}), where we refer to \cite{Bonito2016} and \cite{Hohage2008}, and the references therein, for concrete application examples. In this situation $u_f:\domain\to\IR$ describes a steady-state heat flow, i.e. the unique solution to the elliptic PDE
\begin{align}\label{eq:IntroDarcy}
	\begin{dcases}
		\hfill\divergence(f\cdot\nabla u) &= g, \text{ on $\domain$},\\
		\hfill u\hspace{1mm} &= 0, \text{ on $\partial\domain$},
	\end{dcases}
\end{align}
with $\nabla$ denoting the gradient of a function and $\divergence$ the divergence operator applied to a vector field, respectively. Here, $g>0$ is a \textit{heat source function} on $\domain$, which is known in advance, and $f:\domain\to(0,\infty)$ describes an unknown \textit{conductivity} (or \textit{diffusion}) coefficient function. In real-world observation schemes, idealized observations that exactly match the mathematical formalism of a PDE are typically unavailable, as measurement errors introduce statistical noise, motivating the observation scheme \cref{eq:IntroData} in this work. More precisely, the data consists of two parts - the exact solution $u_f$ of the Darcy problem \cref{eq:IntroDarcy}, corrupted by the second part - the Gaussian noise.
Nonetheless, the goal remains the same: to reconstruct the underlying unknown coefficient function $f$ that fits the pattern of \cref{eq:IntroDarcy}.\\

In the present work, we follow the approaches in \cite{MonardNicklPaternain_2020} and \cite{NicklGeerWang2020} and embed the regression model \cref{eq:IntroData} into a more generic setting. Specifically, given a probability space $(\calZ,\scrZ,\zeta)$ and a possibly nonlinear operator (forward map) $\scrG:\IL^2(\domain)\to\IL_\zeta^2(\calZ)$, we assume to have access to observations $(Y_i,Z_i)_{i=1}^N$, which arise from the equation
\begin{equation}\label{eq:IntrogenData}
	Y_i = \scrG(\theta)(Z_i)+\sigma\varepsilon_i,\, i=1,\dots,N,
\end{equation}
where, in this situation, $Z_1,\dots,Z_N$ are i.i.d. drawn from $\zeta$. The goal is now to provide an estimation strategy for $\theta$, which defines an object in a Sobolev space $H^\alpha(\domain)$. Note that \cref{eq:IntroData} corresponds to the special case where $\scrG$ is the nonlinear solution map $f\mapsto u_f$ of the Darcy problem. The reconstruction of parameters $\theta$ following an indirect and noisy observation scheme as in \cref{eq:IntrogenData} has a long mathematical history. For an extensive overview of models with deterministic noise and a linear relationship $\theta\mapsto \scrG(\theta)$, we refer to \cite{EnglHankeNeubauer2000}.  In the last twenty years a wide field has been established to also consider nonlinear and statistical inverse problems. We refer to \cite{Benning_2018, Kaltenbacher_2008,Dashti_2013,Dashti_2015,Stuart_2010,Nickl_2023,Arridge_2019,Bissantz2007}, and \cite{Kaipio_2005}, as well as to the references therein. For a novel proof in the Bayesian approach within this PDE-constrained regression setting, we refer to \cite{MonardNicklPaternain_2020} and \cite{Monard2021}, where statistical guarantees, especially for the posterior mean, are derived under Gaussian process priors on the parameter of interest $\theta$ and analytical assumptions on the forward map $\scrG$. The applicability to PDE models is illustrated not only for the inversion of non-Abelian X-ray transforms but also for a model as in \cref{eq:IntroData} with $u_f:\domain\to\IR$ being the unique solution of the stationary Schrödinger equation. \cite{Giordano_2020} adopt these ideas for the Darcy problem \cref{eq:IntroDarcy} and provide results on posterior convergence. General assumptions to show posterior convergence are formulated in \cite{Nickl_2023}. In \cite{NicklGeerWang2020}, the authors propose a (generalized) penalized least squares estimator $\hat{\theta}$ for $\theta$ in a Gaussian white noise model, which can be interpreted as minimizer of the quantity
 \begin{equation*}
 	H^\alpha(\domain)\ni \theta\mapsto \norm{Y-\scrG(\theta)}^2+\rate^2\norm{\theta}^2_{H^\alpha(\domain)}\in\pRZz.
 \end{equation*}
 They make use of M-estimation techniques to derive statistical guarantees for its prediction loss. 
 The results are then applied for the special case of considering \cref{eq:IntroData} with $u_f:\domain\to\IR$ being the unique solution of the Darcy problem \cref{eq:IntroDarcy} and also for the time-independent Schrödinger equation. In particular, they provide upper bounds for the prediction loss and show that these are in some situations minimax-optimal and coincide with rates known for linear operators. \magenta{While the white noise model provides a mathematically convenient formulation for theoretical analysis, particularly in deriving statistical guarantees, it is often infeasible in real-world applications. More realistic observation schemes involve discrete measurements at deterministic design points, which can be chosen, for instance, to optimally cover the domain under certain regularity assumptions (see \cite{Daon} and \cite{Teckentrup} for instance). In our framework, we instead consider a random design regression model, which is analytically convenient and aligns well with established approaches in nonparametric statistical theory:} For \textit{random} design regression models described by \cref{eq:IntrogenData}, \cite{BohrNickl2023} and \cite{NicklWang2023} provide statistical and computational guarantees for \magenta{a Tikhonov-regularized estimator, defined as the $D$-dimensional minimizer of the quantity}
\begin{equation}\label{eq:IntroMAP}
	H^\alpha(\domain)\ni\theta\mapsto \frac{1}{N}\sum_{i=1}^N |Y_i-\scrG(\theta)(Z_i)|^2 + \rate^2\norm{\theta}_{H^\alpha(\domain)}^2\in\pRZ.
\end{equation}
Under fairly general assumptions on the forward map $\scrG$, the authors demonstrate the consistency and stability of $\hat{\theta}$, provided that $\theta$ can be sufficiently well approximated by its orthogonal projection onto a finite-dimensional subspace. The key of the proof is an empirical process argument, which implicitly requires the forward map $\scrG$ to satisfy the 'global' stability estimates
\begin{equation}\label{eq:INtrostability}
		\norm{\scrG(\theta)-\scrG(\theta')}_{\IL^2} \lesssim \norm{\theta-\theta'}_{\IL^2}\text{ and }\norm{\scrG(\theta)-\scrG(\theta')}_\infty \lesssim \norm{\theta-\theta'}_\infty
\end{equation}
uniformly over the whole parameter space. While the solution map of the time-independent Schrödinger equation meets these requirements, there is no known analogous result for the second condition for the solution map of the Darcy problem, which renders the theory of \cite{BohrNickl2023} and \cite{NicklWang2023} inapplicable. In addition, the first requirement leads to slower convergence rates, which we improve by taking the ill-posedness of the forward map $\scrG$ into account (see \cite{NicklGeerWang2020} for corresponding results in the white noise model). \magenta{Furthermore, the authors discuss the relationship between $\hat{\theta}$ and the Bayesian approach for inverse problems, where $\hat{\theta}$ can also be interpreted as maximizer of the likelihood $\pi(\cdot|(Y_i,Z_i)_{i=1}^N)$ of the posterior distribution $\Pi(\cdot|(Y_i,Z_i)_{i=1}^N)$, under a Gaussian process prior $\Pi(\cdot)$ on $\theta$. This also results in the naming of $\hat{\theta}$ as maximum a posteriori (MAP) estimator.} In section 2.2. of \cite{NicklWang2023}, the authors propose conditions under which the stochastic gradient algorithm
\begin{equation}
	\vartheta_0 := \theta_{\operatorname{init}}, \hspace{1ex} \vartheta_{k+1} = \vartheta_{k} + \rho\nabla\log \tilde{\pi}(\vartheta_{k}|(Y_i,Z_i)_{i=1}^N)
\end{equation}
converges geometrically to the minimizer $\hat{\theta}$ of \cref{eq:IntroMAP} (and hence to the ground truth $\theta$ after sufficiently many iterations) in polynomial time,  if the initial value $\theta_{\operatorname{init}}\in\IR^D$ is chosen appropriately. We refer to section 2.2.5. (and particularly Theorem 2.8) in \cite{NicklWang2023}. The results are then applied to $\scrG$ being the solution map $f\mapsto u_f$ of the stationary Schrödinger equation, which meets global stability conditions \cref{eq:INtrostability}. Under more general assumptions, \cite{Nickl_2023} discusses situations in which the analysis of the local concavity of the associated log-likelihood to \cref{eq:IntrogenData} is sufficient to specify a stochastic gradient algorithm, for which $\hat{\theta}$ is computable in polynomial time, if the underlying parameter space is finte-dimensional. \magenta{Building on the proof strategy of \cite{NicklWang2023}, the authors of \cite{Altmeyer_2024} establish nonasymptotic sampling guarantees for a gradient-based MCMC algorithm relative to a surrogate posterior in high-dimensional generalized linear models (see Section 2.4), a framework that also encompasses the Darcy problem discussed in Section 3.4. Hence, we refrain from discussing the computational aspects in more detail in this work due to the similarities with \cite{Nickl_2023,NicklWang2023} and \cite{Altmeyer_2024}}.\\

A different approach of modeling real-world phenomena utilizes the description of dynamical (time-dependent) systems via stochastic partial differential equations (SPDEs), which inherently take into account the noise in the system dynamics. The SPDE framework describes uncertainty by embedding stochastic processes directly into the model equations, thereby offering an alternative for capturing the evolution of systems influenced by random fluctuations. For an overview of recent work of SPDEs we refer to \cite{Altmeyer2020,Altmeyer2021} and \cite{gaudlitz2022estimation} as well as the references therein. We emphasize that SPDEs are not the focus of the present work.\\
In the present paper, we follow the approach of \cite{MonardNicklPaternain_2020,Monard2021} and \cite{NicklGeerWang2020} and consider the general nonlinear Gaussian regression model with random design (as \cref{eq:IntrogenData}) and define a maximum a posteriori estimator, similar to \cref{eq:IntroMAP}. We introduce conditions on the forward map $\scrG$, generalizing \cref{eq:INtrostability} by replacing the right-hand sides in both stability estimates with $(H^\kappa)^*$- and $C^\eta$-norms, respectively, while also allowing polynomial growth. We present a general concentration result, where we compare the prediction loss $\frakd^2(\hat{\theta},\theta)$ with some proxy $\frakd^2(\thetatilde,\theta)$ (interpretable as approximation error) and obtain a corresponding consistency result for the MAP estimate, whenever the proxy approximates the ground truth sufficiently well. This result extends the findings of \cite{BohrNickl2023} and \cite{NicklWang2023}, as we allow a generic reference function $\thetatilde$, rather than only a finite-dimensional approximation of the ground truth and provide faster convergence rates for the squared prediction loss, which are now upper-bounded by $N^{-\frac{2(\alpha+\kappa)}{2(\alpha+\kappa)+d}}$,
where $\kappa\in\pRZ$ indicates the degree of ill-posedness of the underlying inverse problem. Under an additional inverse stability assumption on the forward map $\scrG$, an analogous result for the estimation loss of $\hat{\theta}$ is established. For the sake of completeness, we provide an existence result for the MAP-estimator in this random design regression model by exploiting the \textit{direct method in the calculus of variations}, where we adapt a corresponding proof of \cite{NicklGeerWang2020} for the white noise model. The general theory is applied on the Darcy problem. Specifically, we demonstrate that the corresponding solution map satisfies the new stability assumption with $\eta = 1$. In this particular situation, we provide consistency for the MAP estimator, as well as corresponding convergence rate for prediction and estimation errors, which fills a gap in existing theory.\\
\phantom{x}\\

\textbf{Organisation of the paper:} In \cref{sec:model} we introduce the most essential notation and the statistical model in detail. Further, we introduce and discuss the assumption made on the forward map $\scrG$. We define the MAP-estimator, for which we provide an existence result. The main results - consisting of the concentration inequality, an upper bound for the prediction error, as well as consistency for the MAP-estimator - are presented. In \cref{sec:PDE} we apply the main result to the special case of the Darcy problem. We provide analogously consistency for the MAP-estimator, as well as upper bounds for the prediction and estimation risk. In \cref{sec:app}, we provide detailed proofs for the results in the main part and further helping results.\\

%% file: content/model/model.tex
\section{A general statistical inverse Problem}\label{sec:model}

\begin{notation}[Preliminary]
	Throughout the whole paper, we use the following notation.
	\begin{enumerate}[label = \roman*)]
		\item For $\IK\in\{\IN,\IR\}$ and $a\in\IR$ set $\IK_{\geq a} :=\IK\cap[a,\infty)$ and $\IK_{> a}:=\IK\cap (a,\infty)$. Analogously, $\IK_{\leq a}$ and $\IK_{< a}$ are defined. We set $\IN_0:=\IN\cup\{0\}$.
		\item For all $a,b\in\IZ$ with $a\leq b$ set $\nnset{a}{b}:=[a,b]\cap\IZ$, and $\nset{n}:=\nnset{1}{n}$ for any $n\in\IN$.
		\item For all $a,b\in\IR$ define $a\lor b :=\max\{a,b\}$ and $a\land b:=\min\{a,b\}$. 
		\item  For a topological space $(S,\topo)$ we denote by $\scrB_S$ the Borel-$\sigma$-field of $S$ generated by the topology $\topo$. We write $S^d:=\bigtimes_{i\in\nset{d}} S$ for some $d\in\IN$ and denote by $\scrB^{\otimes d}_S:= \bigotimes_{i\in\nset{d}}\scrB_S$ the product $\sigma$-algebra on $S^d$. The topological dual of $S$ is denoted by $S^*$ and consists of all linear and bounded functionals $L:S\to\IR$. If $S$ is equipped with a metric $d_S$, we write $\IB_S(s_o,r):=\{s\in S:\, d_S(s,s_o)\leq r\}$ for the closed ball of radius $r\in\pRZ$ around $s_o\in S$ in $S$.
		\item  For two normed spaces $(S_1,\norm{\cdot}_{S_1})$ and  $(S_2,\norm{\cdot}_{S_2})$, we write $S_1\hookrightarrow S_2$, if $S_1$ is continuously embedded into $S_2$.
		\item For $d\in\IN$ denote by $\norm{\cdot}_{\IR^d}$ the usual Euclidean norm on $\IR^d$. 
		\item Universal constants are denoted by $\const$. If a constant depends on a family of objects $A$, we write $\const(A)$. Constants arising from model assumptions are denoted by $\Const$ with some further specification. In each statement and also in the corresponding proofs, we will number constants consecutively, always starting anew. If not further mentioned, the value of a constant can change from line to line.
	\end{enumerate}
\end{notation}

\input{content/model/funcsec.tex}
\input{content/model/gaussianregmodel.tex}

\subsection{Maximum a posteriori estimator}

 Given the data set $D_N = (Y_i,Z_i)_{i\in\nset{N}}\in(V\times\calZ)^N$ arising from \cref{eq:Data}, fix some $\alpha,\rate\in\pRZz$. We define a penalized \textit{least squares functional} by 
\begin{equation}\label{eq:TikhonovFunc}
	\scrJ_{\rate,N}:	\Thetatilde \to \IR_{\leq 0}\cup\{-\infty\},\text{ }\scrJ_{\rate,N}[\theta]:=-\frac{1}{2\sigma^2N}\sum_{i\in\nset{N}}\Vnorm{Y_i-\scrG(\theta)(Z_i)}^2- \frac{\rate^2}{2}\norm{\theta}_{H^\alpha(\domain,W)}^2
\end{equation}
and set
\begin{equation}\label{eq:MAP}
	\thetaMAP :\in \argmax_{\theta\in\Thetatilde}\left\{\scrJ_{\rate,N}[\theta]\right\}.
\end{equation}

\begin{remark}[Maximum a posteriori estimator]
	\mbox{}
	\begin{enumerate}[label = \roman*)]
		\item Up to now, we do not have formulated any specific conditions on $\Thetatilde$, except of $\Thetatilde\subseteq\IL^2(\domain,W)$. In \cref{propo:existence}, \magenta{we need to restrict ourselves to suitable subsets $\Theta\subseteq\Thetatilde\cap\scrH$ (with $\scrH$ being either $H^\alpha(\domain,W)$ or $H^\alpha_c(\domain,W)$)} in order to provide existence of $\thetaMAP$ as a maximizer of the optimization problem \cref{eq:MAP} \magenta{over $\Theta$}. However, the results of the statistical analysis hold true for \textit{any} maximizer $\thetaMAP$ of $\scrJ_{\rate,N}$.
		From the application perspective, the definition of $\thetaMAP$ seems unfeasible as the optimization problem is formulated over a possibly infinite-dimensional space. The theory shown in \cite{Nickl_2023} and \cite{NicklWang2023} provide general conditions and algorithms, under which $\thetaMAP$ can be computed in polynomial time, if $\Theta$ is finite dimensional. 
		
		\item \magenta{ The estimator $\thetaMAP$ introduced in \cref{eq:MAP} admits two different but related interpretations. On the one hand, it is closely connected to the classical theory of Tikhonov-regularized estimators, as studied in \cite{EnglHankeNeubauer2000}, where the regularization of ill-posed deterministic inverse problems is achieved by penalizing solutions with respect to a suitable Sobolev norm. Specifically, $\thetaMAP$ minimizes the empirical squared loss subject to a regularization term given by the $H^\alpha(\domain, W)$-norm, thereby enforcing smoothness to exclude overly complex solutions.
		On the other hand, $\thetaMAP$ corresponds to the maximum a posteriori (MAP) estimator in a Bayesian inverse problem framework: Let $\Pi$ be a Gaussian process prior on the parameter $\theta$ with reproducing kernel Hilbert space (RKHS) $\scrH$ and corresponding prior density $\pi(\theta)\propto\exp(-\frac{\overline{\rate}^2}{2}\norm{\theta}_{H^\alpha(\domain,W)}^2)$ with $\overline{\rate}>0$. Under the Gaussian likelihood function of the model \cref{eq:Data}, any maximizer $\thetaMAP$ of $\scrJ_{\rate,N}$ over $\scrH$ with $\rate = \frac{1}{\sqrt{N}}\overline{\rate}$ has a formal interpretation as maximizer of the log-likelihood $\log\pi(\theta|D_N)$ of the resulting posterior distributions $\Pi(\cdot| D_N)$. This interpretation has been formalized, for instance, in Remark 4 of \cite{NicklGeerWang2020}. We refer also to \cite{BohrNickl2023,Dashti_2013}, and \cite{Helin_2015} for additional background on MAP estimation in Gaussian Bayesian inverse problems. Throughout this work, we refer to $\thetaMAP$ as the maximum a posteriori (MAP) estimator.}
	\end{enumerate}
\end{remark}

\subsection{A convergence rate result}

The main results of this work - \cref{propo:existence}, \cref{thm:maintheorem} and \cref{co:MAPConsistency} - are presented in the following. Firstly, we provide existence of $\thetaMAP$ as a solution to the optimization problem \cref{eq:MAP} over a suitable subset $\Theta\subseteq\Thetatilde\cap H^\alpha(\domain,W)$. Moreover, we are showing a concentration result controlling the \textit{prediction error} and the regularity of $\thetaMAP$ simultanously under $\IP_{\theta_\frako}^N$-probability. Therefore, we define for every $\theta_1\in\Thetatilde\cap H^\alpha(\domain,W)$, $\theta_2\in\Thetatilde$, and $\rate\in\pRZz$ the functional
\begin{equation*}
	\dist_\rate^2(\theta_1,\theta_2) := \norm{\scrG(\theta_1)-\scrG(\theta_2)}_{\IL^2_\zeta(\calZ,V)}^2 + \rate^2\norm{\theta_1}_{H^\alpha(\domain,W)}^2.
\end{equation*}
We start with analytical assumption on the forward map $\scrG$.
\begin{assumption}[Forward Map I]\label{ass:ForwardMapI}
	\mbox{}
	\begin{enumerate}[label = {C}\arabic*)]
		\item\label{item:first} $\IL^2$-Lipschitz continuity: Let $\alpha,\gamma_1,\kappa\in\pRZz$. Set $\scrH =H^\alpha_c(\domain,W),\text{ if $\kappa\geq \frac{1}{2}$ }$ and $\scrH =	H^\alpha(\domain,W),\text{ if $\kappa < \frac{1}{2}$} $. The map $\scrG$ is $(\alpha,\gamma_1,\kappa)$-regular in $\IL^2$-norm, if there exists a constant $\Const_{\Lip,2}\in\pRZ$, such that for all $\theta_1,\theta_2\in\Thetatilde\cap\scrH$
		\begin{equation*}%\label{eq:AssFMi}
			\norm{\scrG(\theta_1)-\scrG(\theta_2)}_{\IL^2_\zeta(\calZ,V)} \leq \Const_{\Lip,2}\left(1+\norm{\theta_1}_{H^\alpha(\domain,W)}^{\gamma_1}\hspace{-1.5mm}\lor\norm{\theta_2}_{H^\alpha(\domain,W)}^{\gamma_1}\right) \norm{\theta_1-\theta_2}_{(H^\kappa(\domain,W))^*}.
		\end{equation*}
		\item Uniform boundedness: There exists a constant $\Const_{\scrG,\operatorname{B}}\in\IR_{\geq 1}$, such that 
			\begin{equation*}%\label{eq:UniformBoundesG}
				\sup_{\theta\in\Thetatilde,z\in\calZ}\left\{\Vnorm{\scrG(\theta)(z)}\right\} \leq \Const_{\scrG,\operatorname{B}}.
			\end{equation*}
		\item$\IL^\infty$-Lipschitz continuity:  Let $\alpha,\gamma_2,\eta_1,\eta_2\in\pRZz$, such that $\alpha>\eta_1\lor \eta_2+\frac{d}{2}$. The map $\scrG$ is $(\alpha,\gamma_2,\eta_1,\eta_2)$-regular in $\IL^\infty$-norm, if there exists a constant $\Const_{\Lip,\infty}\in\pRZ$, such that for all $\theta_1,\theta_2\in \Thetatilde\cap\scrH$
		\begin{equation*}
			\norm{\scrG(\theta_1)-\scrG(\theta_2)}_\infty \leq \Const_{\operatorname{Lip},\infty}\left(1+\norm{\theta_1}_{C^{\eta_1}(\domain,W)}^{\gamma_2}\lor\norm{\theta_2}_{C^{\eta_1}(\domain,W)}^{\gamma_2}\right)\norm{\theta_1-\theta_2}_{C^{\eta_2}(\domain,W)}.
		\end{equation*} 
	\end{enumerate}
\end{assumption}

\begin{remark}[$\IL^2$- and $\IL^\infty$-Lipschitz continuity]\label{rem:C3}
	\magenta{\item The assumptions C1)--C3) imposed on the forward map $\scrG$ are abstract but serve as minimal requirements for establishing statistical guarantees for the MAP estimator $\thetaMAP$. Condition C1) quantifies the degree of ill-posedness through the smoothing parameter $\kappa$, which also directly influences the convergence rates obtained in \cref{thm:maintheorem}. In contrast, C3) is of a more technical nature: it ensures that entropy numbers of images $\{\scrG(\theta):\theta\in H^\alpha(\domain,W,r)\}$, $r\in\IR_{>0}$, of Sobolev balls under $\scrG$ are sufficiently controlled and is particularly tailored to the PDE-constrained setting considered here, such as the Darcy flow. For this purpose, we weaken the $\IL^\infty$-stability condition
		\begin{equation}\label{eq:rem:c3}
		\norm{\scrG(\theta) - \scrG(\theta')}_\infty \lesssim \norm{\theta - \theta'}_\infty,
		\end{equation}
		as used in \cite{BohrNickl2023}, by replacing the uniform norm on the right with the stronger $C^{\eta_2}$-norm and additionally allowing polynomial growth of order $\gamma_2$. The parameters $\eta_1$ and $\eta_2$ do not appear in the convergence rates obtained in \cref{thm:maintheorem} explicitly, but the condition $\alpha > \eta_1 \lor \eta_2 + \frac{d}{2}$ ensures that the relevant embeddings $H^\alpha(\domain, W) \hookrightarrow C^{\eta_i}(\domain, W)$, $i \in \{1,2\}$, hold, so that all norms are well-defined. It is also not surprising that stronger assumptions are needed to establish consistency of the MAP estimator than to prove posterior contraction in the Bayesian setting, such as in \cite{Giordano_2020} for the Darcy problem. This phenomenon is well known in classical statistics (see, e.g., \cite{LeCam1986}) and becomes evident in the general nonlinear random regression framework as well: while \cite{Nickl_2023} derives posterior contraction under the comparatively mild conditions C1) and C2) (with $\gamma_1 = 0$), the analysis in \cite{BohrNickl2023} requires the uniform stability condition \cref{eq:rem:c3} to prove consistency of the MAP estimator. In our setting, the condition \cref{eq:rem:c3} correspond to the special case $\eta_1 = \eta_2 = \gamma_2 = 0$, which is insufficient for the Darcy problem. As we demonstrate in \cref{sec:PDE} and \cref{sec:app}, the associated solution map fulfills C3) with $\eta_1 = \eta_2 = 1$ and $\gamma_2 = 4$.
	}
\end{remark}

The first main results address the open question of existence of $\thetaMAP$ as the maximizer of the penalized Tikhonov-regularized functional \cref{eq:TikhonovFunc}. For the white noise model, the authors in \cite{NicklGeerWang2020} provide a proof using the white noise formulation of the regression model. We adapt the corresponding proof to the random design regression model, also exploiting  the \textit{Direct Method of the Calculus of Variations}, making reference to \cite{Dacorogna2008}. The proof can be found in \cref{sec:app}.
\begin{proposition}[Existence of $\thetaMAP$]\label{propo:existence}
	Suppose that $\scrG:\Thetatilde\to\IL^2_\zeta(\calZ,W)$ satisfies the conditions C1) - C3) imposed in \cref{ass:ForwardMapI} for some $\alpha>\eta_1\lor\eta_2 +\frac{d}{2}$. Given the data $D_N\sim\IP_{\theta_\frako}^N$ for some fixed $\theta_\frako\in\Thetatilde$, let $\Theta\subseteq\tilde{\Theta}\cap\scrH$ be closed in the weak topology of the Hilbert space $\scrH$. Then for all $\rate\in\pRZ$, almost surely under $\IP_{\theta_\frako}^N$, there exists a maximizer $\thetaMAP$ of $\scrJ_{\rate,N}$ over $\Theta$, satisfying 
	\begin{equation*}
		\sup_{\theta\in\Theta}\left\{\scrJ_{\rate,N}[\theta]\right\} = \scrJ_{\rate,N}[\thetaMAP].
	\end{equation*}
\end{proposition}

Let $R,\rate\in\pRZ$ be fixed. To shorten notation, we introduce the event
\begin{equation*}
	\Xi_{\rate,R} := \left\{\dist_\rate^2(\thetaMAP,\theta_{\frako}) \geq \const_1\cdot(\dist_\rate^2(\tilde{\theta},\theta_\frako)+R^2)\right\}
\end{equation*}
for any maximizer $\thetaMAP$ of $\scrJ_{\rate,N}$ over $\Theta$, some reference parameter $\thetatilde\in \Theta$ (with $\Theta\subseteq \Thetatilde\cap\scrH$) and some $c_1\in\IR_{\geq 1}$. To shorten notation, we define the family of symbols
\begin{equation*}
	\Ptt:=\left\{\alpha,\eta_1,\eta_2, \kappa,\gamma_1,\gamma_2,\Const_{\Lip,2},\Const_{\Lip,\infty},\Const_{\scrG,\operatorname{B}},d,d_W,d_V,\sigma\right\}
\end{equation*}
and write $\const(\texttt{P})$ for a constant to emphasize the dependency.

\begin{theorem}[Concentration inequality]\label{thm:maintheorem}
	Suppose that $\scrG:\Thetatilde\to\IL^2_\zeta(\calZ,V)$ satisfies C1)-C3) in \cref{ass:ForwardMapI} with some $\alpha>\max\left\{\eta_1\lor\eta_2+d\lor\frac{d}{2}(1+\gamma_2),\frac{d}{2}\gamma_1-\kappa\right\}$. Let $D_N\sim\IP_{\theta_\frako}^N$ from \cref{eq:Data} with some fixed $\theta_\frako\in\Thetatilde$. Then the following hold.
	\begin{enumerate}[label = \roman*)]
		\item Let $\Theta\subseteq \Thetatilde\cap\scrH$. Then for any $\bar{\const}\in\IR_{\geq 1}$, we can choose $\const_1 = \const_1(\bar{\const},\Ptt)\in\pRZ$ sufficiently large, such that for all $\rate\in\pRZ$ and $R\in\pRZ$, with $R\geq\rate\geq N^{-\frac{1}{2}}$, satisfying
		\begin{equation}\label{eq:main:ass1}
			R^{-1+\frac{d\gamma_1}{2(\alpha+\kappa)}}	\rate^{-\frac{d(\gamma_1+1)}{2(\alpha+\kappa)}} \leq N^{\frac{1}{2}} \tag{C4}
		\end{equation}
		and 
		\begin{equation}\label{eq:main:ass2}
			R^{-2+\frac{d(\gamma_2 +1)}{\alpha-\eta_2}}	\rate^{-\frac{d(\gamma_2 +1)}{\alpha-\eta_2}} \leq N,\tag{C5}
		\end{equation}
		any maximizer $\thetaMAP$ of $\scrJ_{\rate,N}$ over $\Theta$ and $\thetatilde\in \Theta$ satisfies
		\begin{equation}\label{eq:mainconcentration}
			\IP_{\theta_\frako}^N(\Xi_{\rate,R})  \leq\const_2 \cdot\exp(-\bar{\const} \cdot NR^2)
		\end{equation}
		with $\const_2=\const_2(d_V)\in\pRZ$. 
		
		\item  Let $\Theta\subseteq\Thetatilde\cap\scrH$. Then, there exists a constant $\const_3 = \const_3(\Ptt)\in\pRZ$, such that for all $\rate\in\pRZ$ satisfying 
		\begin{equation}\label{eq:main:ass3}
			\rate^{-1}\lor\rate^{-1-\frac{d}{2(\alpha+\kappa)}} \leq N^{\frac{1}{2}}, \tag{C6}
		\end{equation}
		we have for any maximizer $\thetaMAP$ of $\scrJ_{\rate,N}$ over $\Theta$ 
		\begin{equation*}\label{eq:mainupperbound}
			\IE_{\theta_\frako}^N\left[\dist_\rate^2(\thetaMAP,\theta_{\frako})\right]\leq \const_3\cdot\inf_{\thetatilde\in\Theta}\left\{ \norm{\scrG(\thetatilde)-\scrG(\theta_\frako)}_{\IL^2_\zeta(\calZ,V)}^2+ \rate^2(1+\Vert\thetatilde\Vert_{H^\alpha(\domain)}^2) + \frac{1}{N}\right\}.
		\end{equation*}
	\end{enumerate}
\end{theorem}

\begin{corollary}[Consistency]\label{co:MAPConsistency}
	Under the assumptions of \cref{thm:maintheorem}, let $\theta_\frako \in\Theta$ with $\norm{\theta_\frako}_{H^\alpha(\domain,W)} \leq \operatorname{L}\in\pRZ$. Let $\rate\in\pRZ$ be given by
	\begin{equation*}
		\rate := \rate_N := N^{-\frac{\alpha +\kappa}{2(\alpha+\kappa)+d}}.
	\end{equation*}
	Then, for all $\bar{\const}\in\IR_{\geq 1}$, there exists $\const_{4}=\const_4(\bar{\const},\operatorname{L},\Ptt)\in\pRZ$ large enough, such that any maximizer $\thetaMAP$ of $\scrJ_{\rate,N}$ over $\Theta$ satisfies
	\begin{equation*}
		\IP_{\theta_\frako}^N\left(\dist_\rate^2(\thetaMAP,\theta_{\frako}) \geq \const_4 \rate^2_N\right) \leq \const_5\exp(-\bar{\const}\cdot N\rate_N^2)
	\end{equation*}
	with $\const_5 = \const_5(d_V)\in\pRZ$.
	In particular, there is $\const_6 = c_6(d_V,\operatorname{L})\in\pRZ$, such that
	\begin{equation*}
		\sup_{\theta\in\Theta\cap H^\alpha(\domain,W,\operatorname{L})}\left\{\IE^N_{\theta}\left[\norm{\scrG(\thetaMAP)-\scrG(\theta)}_{\IL^2_\zeta(\calZ,V)}^2\right]\right\}\leq \const_6\cdot N^{-\frac{2(\alpha +\kappa)}{2(\alpha+\kappa)+d}}.
	\end{equation*}
\end{corollary}

\begin{remark}[Consistency and prediction error]
	\mbox{}
	\begin{enumerate}[label = \roman*)]
		\item Note, that \cref{thm:maintheorem} does not require $\theta_\frako\in\Theta$. Therefore, the assumption $\theta_{\frako}\in\Theta$ in \cref{co:MAPConsistency} \magenta{could also be weakened by assuming that $\scrG(\thetatilde)$ approximates $\scrG(\theta_\frako)$ in $\IL^2_\zeta(\calZ,V)$ distance} (see equation 2.23 in \cite{BohrNickl2023} for instance).
		\item Evidently, \cref{co:MAPConsistency} provides the consistency of the MAP estimator $\thetaMAP$. \magenta{The rate of convergence for the prediction loss fits the interpretation of classical inverse problems: The smoothness parameter $\kappa$ of the forward map $\scrG$, quantifying the ill-posedness of the problem, influences directly the rate of convergence. It is known that this rate of convergence is minimax-optimal for $\kappa$-smoothing linear forward maps $\scrG$ in various settings, see, e.g., \cite{Cavalier_2008}. In contrast, $\eta_1$ and $\eta_2$ from condition C3) do not appear explicitely in the rate, but ensure by the condition $\alpha>\eta_1\lor\eta_2+d\lor\frac{d}{2}(1+\gamma_2)$ that the model is sufficiently regular.}
	\end{enumerate}
\end{remark}
Without any further assumptions, there is no hope to obtain an analogous upper bound for the estimation risk as well. However, the next corollary demonstrates stability estimates under an inverse continuity modulus of $\scrG^{-1}$.
\begin{assumption}[Forward Map II]\label{ass:ForwardMapII}
	\mbox{}
	\begin{enumerate}[label = {C}7)]
		\item Inverse continuity modulus: Additionally to \cref{ass:ForwardMapI}, define for any $\delta\in\pRZ$ 
		\begin{equation*}
			\Lambda_\delta:=\left\{(\theta_1,\theta_2)\in\Thetatilde^2:\,\hspace{-1mm}\norm{\theta_1}_{H^\alpha(\domain,W)}+\norm{\theta_2}_{H^\alpha(\domain,W)} \leq \operatorname{M},\, 	\norm{\scrG(\theta_1)-\scrG(\theta_2)}_{\IL^2_\zeta(\calZ,V)} \leq \delta \right\}
		\end{equation*}
		for some universal constant $\operatorname{M}\in\pRZ$. There exists a constant $\Const_{\operatorname{M}}\in\pRZ$ and some $0<\tau \leq 1$, such that for $\delta$ small enough  		
		\begin{equation*}
			\sup\left\{\norm{\theta_1-\theta_2}_{\IL^2(\domain,W)}:\,(\theta,\theta')\in\Lambda_\delta \right\} \leq \Const_{\operatorname{M}}\cdot\delta^\tau.
		\end{equation*}
	\end{enumerate}
\end{assumption}

\begin{corollary}[Stability]\label{cor:stability}
	Under \cref{ass:ForwardMapI} and \cref{ass:ForwardMapII} and under the conditions formulated in \cref{thm:maintheorem}, fory every $\bar{\const}\in\pRZ$, there exists $c_7(\Const_{\operatorname{M}}, \bar{\const},\norm{\theta_{\frako}}_{H^\alpha(\domain,W)},\Ptt)\in\pRZ$ large enough, such that any maximizer $\thetaMAP$ of $\scrJ_{\rate,N}$ over $\Theta$ satisfies
	\begin{equation*}
		\IP_{\theta_\frako}^N\left(\norm{\thetaMAP-\theta_\frako}_{\IL^2(\domain,W)}\leq \const_7\cdot \rate_N^\tau,\,\norm{\thetaMAP}_{H^\alpha(\domain,W)}\leq \const_7\right)\geq 1- \const_8\cdot\exp(-\bar{\const}\cdot N\rate_N^2)
	\end{equation*}
	for some $c_8=c_8(d_V)\in\pRZ$.
\end{corollary}

%% file: content/model/funcsec.tex
\subsection{Hölder and Sobolev spaces}

Throughout let $(\calZ,\scrZ,\zeta)$ be a measure space. For $p\in[1,\infty]$ we define the usual space $\IL_{\zeta}^p(\calZ):=\IL^p(\calZ,\scrZ,\zeta)$ of (equivalence classes of) functions $f:(\calZ,\scrZ) \to (\IR,\scrB_{\IR})$ endowed with the norm $\norm{\cdot}_{\IL^p_\zeta(\calZ)}$. Let $(V,\Vnorm{\cdot})$ be a finite-dimensional normed $\IR$-vector space with dimension $d_V := \operatorname{dim}_\IR(V)$. By identifying $V$ canonically with $\IR^{d_V}$, we define the space of all vector fields $f = (f_1,\dots,f_{d_V}):\calZ\to V$ with the integrable components $f_i\in\IL^p_{\zeta}(\calZ)$, $i\in\nset{d_V}$, as the tensor space $\IL^p_\zeta(\calZ,V):=\bigtimes_{i\in\nset{d_V}} \IL^p_\zeta(\calZ)$. Recalling that $\IL^2_\zeta(\calZ)$ defines a Hilbert space with inner product
\begin{equation*}
	\langle f,g\rangle_{\IL^2_\zeta(\calZ)} := \int_\calZ f(x)g(x)\mathrm{d}\zeta(x)
\end{equation*}
for $f,g\in \IL^2_\zeta(\calZ)$, the space $\IL^2_\zeta(\calZ,V)$ also inherits a Hilbert space structure with an inner product given by 
\begin{equation*}
	\langle \cdot,\cdot\rangle_{\IL^2_\zeta(\calZ,V)} := \sum_{i\in\nset{d_V}}\langle \cdot_i,\cdot_i\rangle_{\IL^2_\zeta(\calZ)},
\end{equation*}
where $ \cdot_i$ indicates the $i$-th component. \magenta{Moreover, for any function $f:\calZ\to V$, we denote by $\norm{f}_\infty :=\sup_{z\in\calZ}\Vnorm{f(z)}$ the uniform norm of $f$ over $\calZ$. We do not highlight the underlying space, if it is clear from the context.} Throughout let $d\in\IN$. We denote by $\domain\subseteq\IR^d$ a bounded domain with smooth boundary $\partial\domain$. For any $s\in\IN$ denote by $C^s(\domain)$ the space of $s$-times continuously differentiable functions $f:\domain\to\IR$ with bounded derivatives. By convention, we denote by $C^0(\domain)$ the space of continuous functions on $\domain$. We extend the definition for any $s\in\pRZ\setminus\IN$ by defining $C^s(\domain)$ as the function space of all maps $f:\domain\to\IR$, such that for all multi-indices $\boldbeta\in\IN_0^d$ with $|\boldbeta|\leq \floor{s}$ - the integer part of $s$ - the derivative $D^{\boldbeta} f$ exists and is $s-\floor{s}$-Hölder continuous. A norm on $C^s(\domain)$ is given by
\begin{equation*}
	\norm{f}_{C^s(\domain)}:=\sum_{\boldbeta\in\IN_0^d:|\boldbeta|\leq\floor{s}}\norm{D^{\boldbeta} f}_\infty + \sum_{\boldbeta\in\IN_0^d:|\boldbeta|=\floor{s}}\sup_{x,y\in\domain:x\neq y}\left\{\frac{|D^{\boldbeta} f(x)-D^{\boldbeta} f(y)|}{\norm{x-y}_{\IR^d}^{s-\floor{s}}}\right\}.
\end{equation*}
A function $f:\domain\to\IR$ is called smooth, if it belongs to $C^\infty(\domain):=\bigcap_{s\in\pRZ}C^s(\domain)$. We further denote by $C^\infty_c(\domain)$ its subspace of all smooth functions having compact support. Further let $\calL^d$ be the usual $d$-dimensional Lebesgue measure on $(\IR^d,\scrB_{\IR}^{\otimes d})$ and denote by $\calL_\domain^d$ its restriction on $\domain$. To shorten notation, set $\IL^2(\domain):=\IL^2_{\calL_\domain^d}(\domain)$ and analogously $\IL^2(\domain,V):=\IL^2_{\calL_\domain^d}(\domain,V)$. For $s\in\IN$, we denote by $H^s(\domain):= H^s(\domain,\IR)$ the usual order-$s$ Sobolev space of weakly differentiable functions $f:\domain\to\IR$ with square integrable weak derivatives up to order $s$. By $\norm{\cdot}_{H^s(\domain)}$, we denote its usual norm. This definition extends to positive $s\in\pRZ\setminus\IN$ by interpolation, see \cite{Triebel1983}. By convention, set $H^0(\domain):=\IL^2(\domain)$.  Further define $H^s_c(\domain):=\overline{C_c^\infty(\domain)}^{\norm{\cdot}_{H^s(\domain)}}$ as well as $\bar{H}^s(\domain)$ to be the space of all functions in $H^s(\domain)$, whose support is $\bar{\domain}$.
If $s \neq k+\frac{1}{2}$ for any $k\in\IN$, we have $	\bar{H}^s(\domain) = H^s_c(\domain)$. Negative-order Sobolev spaces are defined as the topological dual of $\bar{H}^s(\domain)$, i.e., $	H^{-s}(\domain):=H^{-s}(\domain,\IR) :=\left(\bar{H}^s(\domain)\right)^*$. Denoting by $	\operatorname{tr}_s: H^s(\domain) \to \bigtimes_{l=0}^{s-1}H^{s-l-\frac{1}{2}}(\domain)$ the higher-order trace operator for any $s\in\IN$, we have $H^s_0(\domain):=\left\{f\in H^s(\domain):\, \operatorname{tr}_s[f]=0\right\} = H^s_c(\domain)$, where we refer to \cite{Lions1967} for detailed definitions and properties. \magenta{Let $(W,\Wnorm{\cdot})$ be another finite-dimensional and normed $\IR$-vector space with dimension $d_W:=\operatorname{dim}_\IR(W)\in\IN$, such that we canonically identify $W$ with $\IR^{d_W}$.} All of the previous spaces of functions defined on $\domain$ with values in $\IR$ can be extended to spaces of the corresponding vector fields $f:\domain \to W$ in the same manner as before. Indeed, if $F(\domain)$ is one of the above function spaces, we set $F(\domain,W):=\bigtimes_{i\in\nset{d_W}}F(\domain)$. If $F(\domain)$ is equipped with some norm $\norm{\cdot}_{F(\domain)}$, a norm on $F(\domain,W)$ is given by $\norm{\cdot}_{F(\domain,W)}^2 := \sum_{i\in\nset{d_W}}\norm{\cdot_i}_{F(\domain)}^2$. Similarly, if $F(\domain)$ is endowed with a Hilbert space inner product $\langle\cdot,\cdot\rangle_{F(\domain)}$,  $F(\domain,W)$ inherits the Hilbert space structure with the inner product $	\langle \cdot,\cdot\rangle_{F(\domain,W)}:=\sum_{i\in\nset{d_W}}\langle \cdot_i,\cdot_i\rangle_{F(\domain)}$. Further, for simplicity, we write for centered balls with radius $r\in\pRZ$ in function spaces $F(\domain)$ and $F(\domain,W)$, $	F(\domain,r):= \IB_{F(\domain)}(0,r)$ and $F(\domain,W,r):= \IB_{F(\domain,W)}(0,r)$, respectively.

%% file: content/model/gaussianregmodel.tex
\subsection{The measurement model}
 Assuming that $(\calZ,\scrZ,\zeta)$ is a probability space, we are interested in a (possibly nonlinear) \textit{forward map}
\begin{equation*}
	\scrG : \Thetatilde\to\IL^2_{\zeta}(\calZ,V),
\end{equation*}
where $\Thetatilde\subseteq\IL^2(\domain,W)$ denotes a parameter space. We emphasize that in this paper we are only concerned with forward maps $\scrG$ that have for each $\theta\in\Thetatilde$ a unique continuous representative in $\IL^2_{\zeta}(\calZ,V)$, so that point evaluations of the form $\scrG(\theta)(z)$ are well-defined for each $z\in\calZ$. We assume to have access to measurements of a random design Gaussian regression model, i.e. for a given sample size $N\in\IN$ and a fixed variance $\sigma^2\in\pRZ$, we observe $D_N := (Y_i,Z_i)_{i\in\nset{N}}\in(V\times \calZ)^N$ arising from the equation 
\begin{equation}\label{eq:Data}
	Y_i= \scrG(\theta)(Z_i) + \sigma\varepsilon_i,\,i\in\nset{N},\, \theta\in\Thetatilde,
\end{equation}
where  $Z_1,\dots,Z_N$ are i.i.d  copies of $Z$, which is drawn from $\zeta$. The error terms $\varepsilon_1,\dots,\varepsilon_N$ are i.i.d. copies of $\varepsilon\sim \operatorname{N}(0,\operatorname{Id}_V)$ with $\operatorname{Id}_V$ denoting the identity operator on the vector space $V$. The design variable $Z$ and the error $\varepsilon$ are assumed to be independent. \magenta{Note, we assume that the noise variance $\sigma^2$ is known in advance. In settings where $\sigma^2$ is unknown, it could be estimated from the data, for instance via the sample variance. A detailed treatment of a unknown variance setting is beyond the scope of the present work and posponed to future work.} The distribution of a generic observation $D:=(Y,Z)$ arising from \cref{eq:Data} is denoted by $\IP_{\theta}$. The joint product law of $D_N$ on $(V\times \calZ)^N$ is denoted as $\IP_{\theta}^N = \bigotimes_{i\in\nset{N}}\IP_{\theta}$. The corresponding expectation operators are denoted by $\IE_{\theta}$ and $\IE_{\theta}^N$, respectively. Keeping the identification of $V$ with $\IR^{d_V}$ in mind, $\IP_{\theta}$ has a probability density function $p_{\theta}$ with respect to the product measure $\calL^{d_V}\otimes \zeta$, namely
\begin{equation*}
	\forall y\in V, z\in \calZ:\,p_{\theta}(y,z) = (2\pi\sigma^2)^{-\frac{d_V}{2}}\exp\left(-\frac{1}{2\sigma^2}\Vnorm{y-\scrG(\theta)(z)}^2\right).
\end{equation*}
Consequently, $\IP_{\theta}^N$ also has a probability density function with respect to $\bigotimes_{i\in\nset{N}}(\calL^{d_V}\otimes \zeta)$ denoted by $p_{\theta}^N$, which reads for all $(y,z)\in(V\times\calZ)^N$ as
\begin{equation*}
	p_{\theta}^N(y,z) = \prod_{i\in\nset{N}}p_{\theta}(y_i,z_i) = (2\pi\sigma^2)^{-\frac{Nd_V}{2}}\exp\left(-\frac{1}{2\sigma^2}\sum_{i\in\nset{N}}\Vnorm{y_i-\scrG(\theta)(z_i)}^2\right).
\end{equation*}
Given the data vector $D_N$, a shifted log-likelihood function $\ell_N:\Thetatilde\to\IR$ is defined by
\begin{equation*}
	\ell_N(\theta) := \ell_N(\theta|D_N):= \log p_{\theta}^N(D_N) + \frac{Nd_V}{2}\log(2\pi\sigma^2) = -\frac{1}{2\sigma^2}\sum_{i\in\nset{N}}\Vnorm{Y_i-\scrG(\theta)(Z_i)}^2
\end{equation*}
for each $\theta\in\Thetatilde$. Analogously, given a single datum $D =(Y,Z)$, we define
\begin{equation*}
	\ell(\theta) := \ell(\theta|D):= \log p_{\theta}(D) + \frac{d_V}{2}\log(2\pi\sigma^2) = -\frac{1}{2\sigma^2}\Vnorm{Y-\scrG(\theta)(Z)}^2.
\end{equation*}
The results obtained below hold under the frequentist assumption that the data $D_N$ from \cref{eq:Data} is generated from a fixed \textit{ground truth} $\theta_\frako\in\Thetatilde$ including the law $\IP_{\theta_\frako}^N$.

%% file: content/pde/pde.tex
\section{Application: PDE constrained regression models}\label{sec:PDE}
	
	We apply the results from \cref{sec:model} on the inverse problems described by the PDE \cref{eq:IntroDarcy} and derive corresponding convergence rates for the parameter $f$. We further discuss the applicability of the setup for other PDE models.
	Throughout this section, we are interested in estimating the unknown scalar function $f:\domain\to\IR$ and corresponding solutions $u = u_f:\domain\to\IR$ of the Darcy problem. For any $d\in\IN$, integer $\alpha>\frac{d}{2}$, a fixed constant $f_{\min}\in\pRZ$, and denoting by $\partial_n$ the outward pointing normal derivative, define the parameter space  
	\begin{equation}\label{eq:PDEpspace}
		\calF:=\calF_{\alpha,f_{\min}}:=\left\{f\in H^\alpha(\domain):\,f>f_{\min}\text{ on $\domain$}, f|_{\partial\domain}\equiv 1,\partial^j_{n}f \equiv 0,j\in\nset{\alpha-1}\right\},
	\end{equation}
	and its subclass $\calF(\operatorname{L}):=\calF\cap H^\alpha(\domain,\operatorname{L})$, $\operatorname{L}\in\pRZ$. Let $g\in C^\infty(\domain)$ be a known \textit{heat source} function. If $f\in\calF_{\alpha,f_{\min}}$ for some $\alpha>1+\frac{d}{2}$ and $f_{\min}\in\pRZ$, then well-known PDE theory (see \cite{Gilbarg_2001} for instance) provides the existence and uniquess of a classical solution $u\in C^0(\overline{\domain})\cap C^{2+\nu}(\domain)$ of the Dirichlet boundary value problem
		\begin{align}\label{eq:darcy}
			\begin{dcases}
				\hfill\divergence(f(x)\cdot\nabla u(x)) &= g, \text{ on $\domain$},\\
				\hfill u(x)\hspace{1mm} &= 0, \text{ on $\partial\domain$},
			\end{dcases}
		\end{align}
	for some (small) $\nu\in\pRZ$. \magenta{In other words, we can define a non-linear solution map $G^D:\calF\to \IL^2(\domain)$, mapping any $f\in\calF$ to the corresponding solution $u_f = G^D(f)\in\IL^2(\domain)$ of \cref{eq:darcy}.} Following \cref{eq:Data}, we then observe
	\begin{equation}\label{eq:DarcyData}
		Y_i = G^D(f_\frako)(X_i) + \sigma\varepsilon_i, i\in\nset{N},
	\end{equation}
	where $X_1,\dots,X_N$ are drawn from the uniform law on the domain $\domain$, where we assume for simplicity that $\operatorname{vol}(\domain) =1$.
	The inverse problem consists of recovering the ground truth $f_\frako\in\calF$ by using a data set $D_N={(Y_i,X_i)}_{i\in\nset{N}}\in(\IR\times\IR)^N$, whose joint product law is denoted by $\IP_{f_{\frako}}^N$. We apply the results of \cref{sec:model} to a suitable re-parametrization $\Theta$ of $\calF$, which then defines a linear subspace of $H^\alpha(\domain)$. Therefore, we introduce \textit{link functions} $\Psi$.
	\begin{definition}[(Regular) Link Function]\label{def:LinkFunction}
		Let $f_{\min}\in\pRZ$.
		\begin{enumerate}[label = \roman*)]
			\item A function $\Psi:\IR\to\IR_{>f_{\min}}$ is called a \textit{link function}, if $\Psi$ is a smooth, strictly increasing bijective map with $\Psi(0) = 1$ and $\Psi'(x) > 0$ on $\IR$.
			\item A link function $\Psi$ is called \textit{regular}, if and only if all their derivatives are bounded, i.e.
			\begin{equation*}
				\sup_{x\in\IR}\left\{|\Psi^{(k)}(x)|\right\} < \infty,
			\end{equation*}
			for all $k\in\IN$.
		\end{enumerate}
	\end{definition}
	For examples and typical choices of (regular) link functions in practice, we refer to \cite{NicklGeerWang2020} and \cite{NicklWang2023}. In the notation of \cref{sec:model}, throughout this section we set $\calZ = \domain$, $\scrZ :=\scrB_\domain$ and $\zeta = \calL_\domain$ (recalling that $\operatorname{vol}(\domain) = 1$). We set $	\Thetatilde :=\Theta :=\left\{\Psi^{-1}\circ f:\, f\in\calF\right\}= H^\alpha_c(\domain)$ (the \textit{pulled-back}) parameter space, and 
	\begin{equation*}
		\scrG^D : \Theta \to \IL^2(\domain), \, \scrG(\theta):= G^D(\Psi\circ \theta).
	\end{equation*}
	In other words, for each $\theta\in \Theta$, we identify some $f\in \calF$ by $f = \Psi\circ \theta$. To emphasize the identification, we sometimes also write $f= f_\theta$. Given the data set $D_N\sim\IP_{f_{\frako}}^N$, we define a Tikhonov-regularized functional
	\begin{equation}\label{eq:PDETikhonov}
		J_{\rate,N}: \calF\to \IR_{\leq 0}\cup\{-\infty\},\,	f \mapsto -\frac{1}{2\sigma^2N}\sum_{i\in\nset{N}}|Y_i-G^D(f)(X_i)|^2- \frac{\rate^2}{2}\norm{\Psi^{-1}\circ f}_{H^\alpha(\domain)}^2
	\end{equation}
	for any $\rate\in\pRZ$. In the notation of \cref{eq:TikhonovFunc}, we have $J_{\rate,N}[f] = \scrJ_{\rate,N}[\theta]$ for any $f=f_\theta\in\calF$, such that maximizing $J_{\rate,N}$ over $\calF$ is equivalent to maximizing $\scrJ_{\rate,N}$ over $\Theta = H^\alpha_c(\domain)$. Any pairs of maximizers will be denoted by 
	\begin{equation*}
		\fMAP:\in\argmax_{f\in\calF}\left\{J_{\rate,N}[f]\right\} \text{ and } 	\thetaMAP = \Psi^{-1}\circ\fMAP\in\argmax_{\theta\in H^\alpha_c(\domain)}\left\{\scrJ_{\rate,N}[\theta]\right\}.
	\end{equation*}
	 For a link function $\Psi$, define for $f_1,f_2\in\calF$
	\begin{equation*}
		\Dist_\rate^2(f_1,f_2):=\norm{G(f_1)-G(f_2)}_{\IL^2_\zeta(\calZ)}^2+\rate^2\norm{\Psi^{-1}\circ f_1}_{H^\alpha(\domain)}^2.
	\end{equation*}
	Identifying $f_{\theta_1}$ and $f_{\theta_2}$ with $\theta_1,\theta_2\in H^\alpha_c(\domain)$, respectively, we see
	\begin{equation*}
		\Dist_\rate^2(f_{\theta_1},f_{\theta_2}) =	\dist_\rate^2(\theta_1,\theta_2).
	\end{equation*}

\input{content/pde/darcy.tex}

	\subsection*{Comments on further models}
		The general framework presented in \cref{sec:model} does not only apply to the Darcy problem, it promises to consider various PDE-related problems. As already brielfly discussed in \cref{rem:C3}, the results in \cite{BohrNickl2023} concerning the MAP estimate are special cases by setting $\eta_1 = \eta_2=\gamma_2=0$. Hence, the solution map of the time-independent Schrödinger equation satisfies the requirements of \cref{ass:ForwardMapI} and \cref{ass:ForwardMapII} taking the results of \cite{NicklGeerWang2020} also into account. Moreover, very analogous and elementary computations show that the solution map of the heat equation studied in \cite{Kekkonen2022} also satisfies the conditions in \cref{ass:ForwardMapI} and \cref{ass:ForwardMapII} with $\eta_1 = \eta_2=\gamma_2=0$.\\

		In very recent work, \cite{NicklTiti2024} are considering a non-linear Bayesian data assimilation model for the periodic $d=2$-dimensional Navier Stokes equation with unknown initial condition. In this setting $u=u_\theta:\domain\times (0,T)\to\IR^2$ solves the non-linear evolution equation
		\begin{equation*}
			\begin{dcases}
				\frac{\mathrm{d}u}{\mathrm{d}t} + \nu Au + B(u,u) &= f,\\
				\hfill u(0) & = \theta,
			\end{dcases}
		\end{equation*}
		where $\theta$ is the initial condition of the system, $f$ a forcing term and $\nu\in\pRZ$ a viscosity parameter. For this type of equation in two dimensions a complete solution theory already exists. In particular, if the initial condition $\theta$ is contained in some specific function class $\tilde{\calF}$, the solution operator $\tilde{\calF}\ni \theta\mapsto G^{NS}(\theta) = u_f\in\IL^2(\domain\times(0,T))$ is always well-defined.  Based on a data set $D_N\in(\IR^2\times\domain\times(0,T))^N$ arising from \cref{eq:Data}, it seems feasible to apply the theory from \cref{sec:model}. Indeed, using the theory presented in \cite{NicklTiti2024}, elementary computations show that the solution map is $(\alpha,0,0)$-regular in $\IL^2$-norm and uniformly bounded, thus satisfying C1) and C2) from \cref{ass:ForwardMapI}. C7) from \cref{ass:ForwardMapII} can also be shown with $\tau = 1$. If C3) of \cref{ass:ForwardMapI} is additionally satisfied, one would obtain rates of the form $N^{-\frac{\alpha}{\alpha+2}}$ for estimating $\theta$  - both in prediction and estimation error.

%% file: content/pde/darcy.tex
\begin{theorem}[Consistency and stability - Darcy problem]\label{thm:main:darcy}
	Let $\calF$ be as in \cref{eq:PDEpspace} for some $\alpha>1+\frac{5}{2}d$ and $f_{\min}\in\pRZ$. Let $G^{D}(f) = u_f$ denote the unique solution of \cref{eq:darcy} and let $D_N\sim\IP_{f_{\frako}}^N$ from \cref{eq:DarcyData} for some $f_\frako\in\calF$. Let $\Psi$ be a regular link function. Let $J_{\rate_{N},N}$ be as in \cref{eq:PDETikhonov} with
	\begin{equation*}
		\rate_{N} := N^{-\frac{\alpha+1}{2(\alpha+1)+d}}.
	\end{equation*}
	Then the following hold.
	\begin{enumerate}[label = \roman*)]
		\item For each $f_\frako\in\calF$ and, almost surely under $\IP_{f_\frako}^N$, there exists a maximizer $\fMAP\in\calF$ of $J_{\rate_{N},N}$  over $\calF$.
		\item Let $\operatorname{L}\in\pRZ$. Then there exist constants $c_1 = c_1(d,\sigma,\domain,f_{\min},g,\Psi,\alpha,\operatorname{L}),c_2,c_3\in\pRZ$, such that for any maximizer $\fMAP$, we have for all $\operatorname{M}\geq c_1$
				\begin{equation}
					\sup_{f\in\calF(\operatorname{L})}\left\{\IP_{f}^N\left(\Dist_\rate^2(\fMAP,f)\geq \operatorname{M}\cdot	\rate_N^2 \right)\right\} \leq \const_2 \exp\left(-\const_3 \cdot\operatorname{M}^2N	\rate_N ^2\right).
		    	\end{equation}   
		  \item Assume additionally that $g\geq g_{\min}\in\pRZ$. For every $\bar{\const}\in\pRZ$, there exists a constant $c_2 = c_2(\alpha,d,\domain,f_{\min},g_{\min},f_{\frako},\sigma)$ large enough and $c_3\in\pRZ$, such that
		  \begin{multline*}
		  		\IP_{f_{\frako}}^N\left(\norm{\fMAP-f_\frako}_{\IL^2(\domain)} \leq c_2\cdot\rate_N^{\frac{\alpha-1}{\alpha+1}},\,\norm{\Psi^{-1}\circ\fMAP}_{H^\alpha(\domain)}\leq\const_2\right)\\\geq 1-\const_3\cdot\exp(-c_1\cdot N\rate_N^2).
		  \end{multline*}
	\end{enumerate}
\end{theorem}

\begin{theorem}[Prediction and estimation risk]\label{thm:main2:darcy}
	Under the conditions formulated in \cref{thm:main:darcy}, we have the following results.
	\begin{enumerate}[label = \roman*)]
		\item (Prediction Error) For every $\operatorname{L}\in\pRZ$, $\beta\in[0,\alpha+1]$, there exists a finite constant $c_4 =c_4(\alpha,\beta,d,\domain,f_{\min},L,\norm{g}_{H^{\alpha-1}(\domain)})\in\pRZ$, such that any maximizer $\fMAP\in\calF$ of $J_{\rate_{N},N}$ over $\calF$ satisfies
		\begin{equation*}
			\sup_{f\in\calF(\operatorname{L})}\left\{\IE_{f}^N\left[\norm{G^D(\fMAP)-G^D(f)}_{H^\beta(\domain)}^2\right]\right\}\leq\const_4\cdot N^{-\frac{2(\alpha+1-\beta)}{2(\alpha+1)+d}}.
		\end{equation*}
		\item (Estimation Error) Let $g\geq g_{\min}\in\pRZ$. For every $\operatorname{L}\in\pRZ$, there exists a constant $c_5=c_5(g_{\min},f_{\min},\domain,d,\alpha,L,\norm{g}_{H^{\alpha-1}(\domain)})\in\pRZ$, such that any maximizer $\fMAP\in\calF$ of $J_{\rate_{N},N}$ over $\calF$ satisfies
		\begin{equation*}
			\sup_{f\in\calF(\operatorname{L})}\left\{\IE_{f}^N\left[\norm{\fMAP-f}_{\IL^2(\domain)}^2\right]\right\}\leq\const_5\cdot N^{-\frac{2(\alpha-1)}{2(\alpha+1)+d}}.
		\end{equation*}
	\end{enumerate}
\end{theorem}

\begin{remark}[Prediction and estimation risk]\label{rem:rateoptDarcy}
	\magenta{
	\begin{enumerate}[label = \roman*)]
			\item The convergence rates for the prediction risk presented in \cref{thm:main2:darcy} i) align with the classical interpretation of statistical inverse problems. The smoothness of the parameter of interest $f$ is quantified by $\alpha$, while the ill-posedness of the problem is captured by $\kappa = 1$, reflecting the smoothing property of the forward map $G^D$. Indeed, Lemma 23 in \cite{NicklGeerWang2020} shows that if $f$ is $\alpha$-smooth, then $G^D(f)$ is $\alpha + 1$-smooth. Considering the prediction risk in the $H^\beta$-norm, which corresponds to differentiation of order $\beta$, yields the characteristic numerator $\alpha + 1 - \beta$ in the convergence rate. These rates coincide exactly with those obtained in the white noise model setting in \cite{NicklGeerWang2020} for the Darcy problem. Theorem 10 in \cite{NicklGeerWang2020} in addition shows that for the $H^2$-prediction loss the obtained rate of convergence is minimax-optimal. This optimality also carries over to the random design model considered here, as the proof techniques can be adapted without difficulty. Hence, the results in \cite{NicklGeerWang2020} provides a rigorous benchmark indicating that the convergence rates derived in these settings are indeed statistically sharp. 
			\item The additional assumption $g \geq g_{\min} > 0$ in \cref{thm:main:darcy}~iii) and \cref{thm:main2:darcy}~ii) plays a crucial role in the analysis. It ensures the validity of a stability estimate for the inverse operator $(G^D)^{-1}$, which is central to deriving the convergence rate for the estimation error. Such a condition is also imposed in Lemma~24 of \cite{NicklGeerWang2020}, which serves as a foundation for the current proof strategy. We additionally refer to \cite{Richter_1981} and \cite{Bonito_2017} for a broader discussion on injectivity and global stability estimates for inverse problems of this type.
			\item In contrast to the prediction error, the rates of convergence for the estimation error derived in \cref{thm:main2:darcy} ii) do not match those obtained in \cite{Giordano_2020} for the posterior contraction, which are generally slower. This discrepancy is partly due to the lower posterior regularity compared to the MAP estimator, resulting in suboptimal rates for the Bayesian approach. However, the optimality of the estimation error rates for $\fMAP$ established here remains an open question. Classical nonparametric theory would suggest a faster rate, such as having exponent $\frac{2\alpha}{2(\alpha+1)+d}$. It is unclear to the author whether this gap is a consequence of the estimation methodology and the proof techniques employed—particularly those relying on inverse stability estimates and associated PDE analytic arguments—or whether it reflects an intrinsic structural limitation of the Darcy problem itself. Further investigation is needed to clarify this issue.
	\end{enumerate}}
\end{remark}

The proofs of \cref{thm:main:darcy} and \cref{thm:main2:darcy} are given in detail in \cref{sec:app}. They are based on verifying that the forward map $\scrG^D$ satisfies C1)-C4) formulated in \cref{ass:ForwardMapI} and \cref{ass:ForwardMapII} for suitable parameters $\alpha,\kappa,\gamma_1,\gamma_2,\eta_1$ and $\eta_2$, with $V=W=\IR$, and an application of the results in \cref{sec:model}. 

%% file: content/appendix/app_map.tex
\subsection{Proofs of \cref{sec:model}}
In this section, we provide detailed proofs of the main results presented in \cref{sec:model}. We start with the proof of \cref{propo:existence} - the existence result for the MAP-estimator. To do so, fix $N\in\IN$ as well as $\rate\in\pRZ$ throughout. The proof follows the approach of \cite{NicklGeerWang2020}, where the existence of $\thetaMAP$ is shown for the White Noise formulation of the model. The authors make use of the \textit{direct method of the calculus of variations} (see \cite{Dacorogna2008} for instance) and split the proof into two parts: Firstly, it is shown that $\thetaMAP$ exists as a maximizer of the functional $\scrJ_{\rate,N}$ over bounded subsets of $\Theta$ (local existence). Secondly, it is shown that this already suffices to conclude the existence of $\thetaMAP$ as a maximizer over the whole parameter space $\Theta$ (global existence). In the present work, we adapt the proof idea to the random design regression setting. We need the functional $-\scrJ_{\rate,N}$ to be lower semicontinuous with respect to the weak topology of Sobolev spaces.\\ Let $\tilde{\scrH}$ be either $H^\alpha(\domain)$ or $H^\alpha_c(\domain)$. Let $\topo_w$ denote the weak topology of $\tilde{\scrH}$, i.e., the coarsest topology with respect to which all bounded linear functionals $L:\tilde{\scrH}\to\IR$ are continuous. With a slight abuse of notation, the corresponding subspace topology on subsets of $\tilde{\scrH}$ is also denoted by $\topo_w$. By standard theory of functional analysis (see for instance \cite{Megginson1998}, Theorem 2.6.23), on any closed ball $\tilde{\scrH}(r):=\{\theta\in \tilde{\scrH}:\,\norm{\theta}_{H^\alpha(\domain)}\leq r\}$ of radius $r\in\pRZ$, the weak topology $\topo_w$ is metrizable by some metric $\tilde{d}$. Defining the product space $\scrH :=\bigtimes_{i\in\nset{d_W}}\tilde{\scrH}$, the weak topology $\topo_w^W$ on $\scrH$ equals to $\bigtimes_{i\in\nset{d_W}}\topo_w$ (see proof of Theorem 3.10 in \cite{Brezis2011}). As above, on each closed ball $\scrH(r):=\{\theta\in \scrH:\,\norm{\theta}_{H^\alpha(\domain,W)}\leq r\}$, the weak topology $\topo_w^W$ is then also metrizable by some metric $d$.\\ 
\phantom{x}\\
In the following we claim that $- \scrJ_{\rate,N}:\Theta\cap\scrH(2^j)\to\IR\cup\{\infty\}$ is lower semicontinous with respect to $\topo_w^W$ for any $j\in\IN$. To this end, we need the following Lemma.

\begin{lemma}[Compact Operators and Weak Convergence]\label{lem:compactoperator}
	Let $(X,\norm{\cdot}_X)$ and $(Y,\norm{\cdot}_Y)$ be two Banach spaces. Let $T:X\to Y$ be a linear and compact operator. Assume that $(x_n)_{n\in\IN}\subseteq X$ is a weakly convergent sequence with limit $x\in X$. Then, $(T[x_n])_{n\in\IN}\subseteq Y$ converges strongly in $Y$, i.e. $\lim_{n\to\infty}\norm{T[x_n]-T[x]}_Y=0$. 
\end{lemma}

\begin{proof}[Proof of \cref{lem:compactoperator}]
	\magenta{
		Let $(x_n)_{n\in\IN}\subseteq X$ be a weakly convergent sequence with limit $x\in X$, i.e., for all $L\in X^*$, it holds $\lim_{n\to\infty} L[x_n] = L[x]$. By an application of the uniform boundedness principle, there exists $B\in\pRZ$, such that $\norm{x_n}_X\leq B$ for all $n\in\IN$. In other words, $(x_n)_{n\in\IN}$ is bounded as well as every sub-sequence of $(x_n)_{n\in\IN}$. By hypothesis on $T$ being compact, the set $T((x_n)_{n\in\IN})$ is relative compact in $Y$, i.e. $\overline{T((x_n)_{n\in\IN})}$ is compact in $Y$. Thus, any subsequence $(T[x_{n_m}])_{m\in\IN}$ of  $(T[x_{n}])_{n\in\IN}$ has a further strongly converging subsequence $(T[x_{n_{m_l}}])_{l\in\IN}$ with limit $y\in Y$, that is $\lim_{l\to\infty}\Vert T[x_{n_{m_l}}]-y\Vert_Y=0$. Further, observe that $(T[x_{n}])_{n\in\IN}$ converges weakly in $Y$, as for any $A\in Y^*$ we have $A\circ T\in X^*$ and thus $\lim_{n\to\infty}A[T[x_{n}]]=\lim_{n\to\infty}(A\circ T)[x_n] = (A\circ T)[x] = A[T[x]]$. In other words, we have shown that $(T[x_n])_{n\in\IN}$ converges weakly to $T[x]$ in $Y$ and $(T[x_{n_{m_l}}])_{l\in\IN}$ converges strongly to $y$ in $Y$. By uniqueness of the weak limit, we deduce $ y = T[x]$ and therefore $T[x_n]$ converges strongly to $T[x]$ in $Y$. }
\end{proof}

\begin{remark}[Embeddings and Weak Convergence]\label{rem:weakconvergence}
	For $d\in\IN$ let $\domain\subseteq\IR^d$ be a bounded domain with smooth boundary. Denoting by $B_{pq}^s(\domain)$ the ususal Besov space with parameters $p,q\in\pRZ\cup\{\infty\}$ and $s\in\IR$, the embedding 
	\begin{equation*}
		B_{p_1q_1}^{s_1}(\domain) \hookrightarrow B_{p_2q_2}^{s_2}(\domain) 
	\end{equation*}
	is compact, if $p_1,p_2,q_1,q_2\in\pRZ\cup\{\infty\}$, $-\infty< s_2 < s_1 < \infty$ and $ s_1-s_2 - d\max\{\frac{1}{p_1}-\frac{1}{p_2},0\}\in\pRZ$ (see \cite{Triebel1996} for definitions and properties). Let $\alpha\in\pRZ$ and $\eta\in \pRZz\setminus\IN$. In particular, by identifying $H^\alpha(\domain) = B_{22}^\alpha(\domain)$, $\IL^2(\domain) = B_{22}^0(\domain)$ and $C^\eta(\domain) = B_{\infty\infty}^\eta(\domain)$, the embeddings $H^\alpha(\domain)\hookrightarrow \IL^2(\domain)$ and $H^{\alpha}(\domain)\hookrightarrow C^\eta(\domain)$ for $\alpha > \eta + \frac{d}{2}$ are compact. Further, following Proposition 4.5 in \cite{Triebel2006}, for any $k\in\IN$ the  embedding $B_{\infty1}^k(\domain)\hookrightarrow C^k(\overline{\domain})$ is continuous. Hence, the embedding $H^\alpha(\domain) \hookrightarrow C^k(\domain)$ is also compact, provided that $\alpha > k+\frac{d}{2}$. If $(W,\Wnorm{\cdot})$ denotes a $d_W$-dimensional $\IR$-vector space, the same results apply, if we replace the function spaces $F(\domain)$ by their tensor products $F(\domain,W)$. \magenta{With $\scrH$ as defined above, let} $(\theta_n)_{n\in\IN}\subseteq \scrH$ be a sequence, such that $\theta_n\to\theta \in \scrH$ with respect to $\topo_w^W$. An application of \cref{lem:compactoperator} then yields that $(\theta_n)_{n\in\IN}$ converges strongly in $\IL^2(\domain,W)$ as well as in $C^\eta(\domain,W)$, if $\alpha > \eta +\frac{d}{2}$ for $\eta\in\pRZz$.
\end{remark}

\begin{lemma}[Continuity of $\scrG$]\label{lem:continuityG}
	For $\gamma_2,\eta_1,\eta_2\in\pRZz$, let $\scrG:\Theta\to\IL^2_\zeta(\calZ,V)$ satisfy C3) of \cref{ass:ForwardMapI}. Let $\Theta_o\subseteq\Theta$ be a bounded subset of $\scrH$ with $\alpha>\eta_1\lor\eta_2+\frac{d}{2}$. Let $(\Omega,\calA, \IP)$ be the underlying probability space of $Z:\Omega \to \calZ$. Then, the map
	\begin{equation*}
		(\Theta_o,d) \ni \theta \mapsto \scrG(\theta)(Z) \in (V,\Vnorm{\cdot})
	\end{equation*}
	is continuous $\IP$-almost surely.
\end{lemma}
\begin{proof}[Proof of \cref{lem:continuityG}]
		Fix some $\omega\in\Omega$ (outside of a $\IP$-null-set). Let $(\theta_n)_{n\in\IN}\subseteq\Theta_o$ be a convergent sequence with weak limit $\theta\in\Theta_o$, i.e., $\lim_{n\to\infty} d(\theta_n,\theta) = 0$. Note that $\norm{\theta_n}_{H^\alpha(\domain,W)}\lor \norm{\theta}_{H^\alpha(\domain,W)}\leq B$ for some $B\in\pRZ$ uniformly over all $n\in\IN$. The claim is proven, if $(\varphi_n(Z(\omega)))_{n\in\IN}\subseteq\IR$, defined by 
		\begin{equation*}
			\varphi_n(Z(\omega)) := \Vnorm{\scrG(\theta_n)(Z(\omega))-\scrG(\theta)(Z(\omega))}
		\end{equation*} 
		is a null-sequence. By hypothesis and exploiting the continuous embedding $H^\alpha(\domain,W)\hookrightarrow C^{\eta_1}(\domain,W)$, there exists a constant $\const\in\pRZ$ (depending on $\alpha,\eta_1,\eta_2,d,\Const_{\Lip,\infty}$), such that
		\begin{align*}
			\varphi_n(Z(\omega)) \leq \norm{\scrG(\theta_n)-\scrG(\theta)}_\infty & \leq \const\cdot(1+B^{\gamma_2})\cdot\norm{\theta_n-\theta}_{C^{\eta_2}(\domain,W)}.
		\end{align*}
		Exploiting \cref{rem:weakconvergence}, we have $\lim_{n\to\infty}\norm{\theta_n-\theta}_{C^{\eta_2}(\domain,W)=0}$, which proves the claim.
\end{proof}

\begin{lemma}[Uniform Continuity]\label{lem:uniformcontinuityG}
	Let $\scrG:\Theta\to\IL^2_{\zeta}(\calZ,V)$ satisfy C2) and C3) of \cref{ass:ForwardMapI}. Let $D_N=((Y_i,Z_i))_{i\in\nset{N}}\sim\IP_{\theta_\frako}^N$ for $N\in\IN$ as in \cref{thm:maintheorem}. Let $r\in\pRZ$. Then, the map 
	\begin{equation*}
		E:(\Theta\cap \scrH(r),d)\to(\IR,|\cdot|),\,\theta \mapsto E[\theta]:=\frac{1}{2N\sigma^2}\sum_{i\in\nset{N}}\Vnorm{Y_i-\scrG(\theta)(Z_i)}^2
	\end{equation*}
	is almost surely uniformly continuous. 
\end{lemma}
\begin{proof}[Proof of \cref{lem:uniformcontinuityG}]
	For any $\delta\in\pRZ$, denote by $M_\delta$ the modulus of continuity of $E$, i.e.,
	\begin{equation*}
		M_\delta :=\sup_{\theta,\thetatilde\in\Theta\cap\scrH(r): d(\theta,\thetatilde)\leq \delta}\left\{|E[\theta]-E[\thetatilde]|\right\},
	\end{equation*}
	which defines a random variable. If $(\Omega,\mathcal{A},\IP)$ denotes the underlying probability space supporting the law of $Z$ and $\varepsilon$, respectivly, we define the measurable event 
	\begin{equation*}
		A:=\left\{\omega\in\Omega:\, \lim_{\delta\to 0} M_\delta(\omega) =0\right\}.
	\end{equation*}
	We aim to show $\IP(A) = 1$. As $M_\delta$ is decreasing in $\delta$, it suffices to verify that $\lim_{\delta\to 0}\IE[M_\delta] = 0$. Indeed, elementary computations show that
	\begin{align*}
		|E[\theta]-E[\thetatilde]| 
		& \leq \frac{2\Const_{\scrG,\operatorname{B}}}{N\sigma^2}\sum_{i\in\nset{N}}\Vnorm{\scrG(\thetatilde)(Z_i)-\scrG(\theta)(Z_i)} + \frac{1}{N\sigma}\sum_{i\in\nset{N}}|\langle\varepsilon_i,\scrG(\thetatilde)(Z_i)-\scrG(\theta)(Z_i)\rangle_V|.
	\end{align*}
	Define for each $\delta\in\pRZ$ and $Z\sim\zeta$ the modulus of continuity of $\scrG$ by
	\begin{equation*}
		R_\delta := \sup_{\theta,\thetatilde\in\Theta\cap\scrH(r): d(\theta,\thetatilde)\leq \delta}\left\{\Vnorm{\scrG(\thetatilde)(Z)-\scrG(\theta)(Z)}\right\},
	\end{equation*}
	which in turn defines also a random variable.
	In \cref{lem:continuityG}, we have shown that $\scrG$ is almost surely continuous as a mapping from $\Theta\cap\scrH(r)$ to $V$. As $\Theta\cap\scrH(r)$ is a compact metric space, $\scrG$ is in fact almost surely uniformly continuous, and hence $\lim_{\delta\to 0} R_\delta = 0 $ almost surely. Moreover, since $Z$ and $\varepsilon\sim N(0,\operatorname{Id}_V)$ are independent by hypothesis, we have 
	\begin{align*}
		\IE\left[\sup_{\theta,\thetatilde\in\Theta\cap\scrH(r): d(\theta,\thetatilde)\leq \delta}\left\{|\langle\varepsilon,\scrG(\thetatilde)(Z)-\scrG(\theta)(Z)\rangle_V |\right\}\right] &\leq \IE[\Vnorm{\varepsilon}R_\delta] = \sqrt{d_V}\cdot \IE[R_\delta].
	\end{align*}
	Now by exploiting C2), $(R_\delta)_{\delta\in\pRZ}$ is uniformly integrable and hence, the last expectation is converging to zero for $\delta\to 0$. Overall, we obtain
	\begin{equation*}
		\IE[M_\delta] \leq \left(\frac{2\Const_{\scrG,\operatorname{B}}}{\sigma^2}+\frac{\sqrt{d_V}}{\sigma}\right)\IE[R_\delta]\overset{\delta\to 0}{\longrightarrow}0,
	\end{equation*}
	which shows the claim.
\end{proof}

\begin{proof}[Proof of \cref{propo:existence}]
	Grant the assumptions of \cref{propo:existence}.	
	We follow along the lines of \cite{NicklGeerWang2020} and split the proof of existence in two parts. 
	\begin{enumerate}
		\item[i)] \textit{Local Existence:} Let $j\in\IN$ be arbitrarily. We show now that  $\scrJ_{\rate,N}$ has a maximizer over $\Theta\cap\scrH(2^j)$. As $\Theta$ is assumed to be weakly closed and $\scrH(2^j)$ is weakly sequentially compact, due to the separability and reflexivity of Sobolev spaces and by an application of the Banach-Alaoglu theorem, any sequence $(\theta_n)_{n\in\IN}\subseteq \Theta\cap\scrH(2^j)$ has a weakly convergent sub-sequence $(\theta_{n_k})_{k\in\IN}$ with weak limit $\theta\in\Theta\cap\scrH(2^j)$. Thus, following the \textit{Direct Method of the Calculus of Variations}, it suffices to show that $-\scrJ_{\rate,N}|_{\Theta\cap\scrH(2^j)}$ is almost surely weakly lower semicontinuous. Rewrite $-\scrJ_{\rate,N}|_{\Theta\cap\scrH(2^j)}$ as 
		\begin{equation}\label{eq:localexistenceI}
			-\scrJ_{\rate,N}|_{\Theta\cap\scrH(2^j)}(\theta) = E_1(\theta)+E_2(\theta) 
		\end{equation}
		with
		\begin{equation*}
			E_1(\theta) := \frac{1}{2N\sigma^2}\sum_{i\in\nset{N}}\Vnorm{Y_i-\scrG(\theta)(Z_i)}^2 \text{ and } 	E_2(\theta):=\frac{\rate^2}{2}\norm{\theta}_{H^\alpha(\domain,W)}^2.
		\end{equation*}
		 $E_2$ is weakly lower semicontinuous, which is a standard result in functional analysis concerning norms on Banach spaces. In \cref{lem:uniformcontinuityG}, we have derived that $E_1$ is almost surely uniformly continuous as a map between the metric spaces $(\Theta\cap \scrH(2^j),d)$ and $(\IR,|\cdot|)$. In particular, $E_1$ is weakly lower semicontinuous. All in all, local existence of $\thetaMAP$ follows immediately by an application of the \textit{Direct Methods of Calculus of Variations} (see \cite{Dacorogna2008}).
		\item[ii)]\textit{Global Existence:} For each $j\in\IN$ and any maximizer $\thetaMAP$ of $\scrJ_{\rate,N}$, we define the event 
		\begin{equation*}
			A_j:=\left\{\omega\in\Omega:\,\thetaMAP(\omega)\notin\Theta\cap\scrH(2^j)\right\}.
		\end{equation*}
		Then, $(A_j)_{j\in\IN}$ satisfies $\lim_{j\to\infty}\IP[A_j] = 0$ according to the concentration inequality in \cref{thm:maintheorem} i). Thus, choosing $j\in\IN$ large enough, we have 
		\begin{equation*}
			\sup_{\theta\in\Theta\cap\scrH(2^j)}\left\{\scrJ_{\rate,N}[\theta]\right\} = 	\sup_{\theta\in\Theta}\left\{\scrJ_{\rate,N}[\theta]\right\}
		\end{equation*}
		with probability as close to one as desired.
	\end{enumerate}
\end{proof}

We are now approaching the proof of \cref{thm:maintheorem}. Therefore, we need some notations and further helping lemma.
\begin{notation}[Empirical Measures]\label{not:empiricalmeasures}
	For $N\in\IN$, let $D_N :=(Y_i,Z_i)_{i\in\nset{N}}\in(V\times \calZ)^N$ be given data.  For any $(y,z)\in V\times\calZ$, denote by $\delta_{(y,z)}(\cdot)$ the Dirac measure on $(V\times\calZ,\scrB_{V\times\calZ})$ in $(y,z)$. The \textit{empirical measure} on $(V\times\calZ,\scrB_{V\times\calZ})$ based on $D_N$ is then for any $A\in \scrB_{V\times\calZ}$ defined by
	\begin{equation*}
		\empM(A) := \frac{1}{N}\sum_{i\in\nset{N}}\delta_{(Y_i,Z_i)}(A).
	\end{equation*}
	Evidently, we have for any measurable function $h:V\times\calZ \to \IR$
	\begin{equation*}
		\int_{V\times\calZ}h(y,z)\mathrm{d}\empM(y,z) = \frac{1}{N}\sum_{i\in\nset{N}}h(Y_i,Z_i).
	\end{equation*}
\end{notation}

The following Lemma provides M-estimation techniques, which are used to prove \cref{thm:maintheorem}. The proof follows along the lines of \cite{Geer_2001},  \cite{NicklGeerWang2020} and \cite{BohrNickl2023}.

\begin{lemma}[M-Estimation Techniques]\label{lem: helpinglemmaI}
	With the notation of \cref{thm:maintheorem} and the lines above it, for any $\theta_{\frako}\in\Thetatilde$, $\thetatilde\in\Theta$, $\rate\in\pRZ$ and any maximizer $\thetaMAP$ of $\scrJ_{\rate,N}$ over $\Theta$, we have
	\begin{equation}\label{eq:MestimationI}
		\frakd_\rate^2(\thetaMAP,\theta_\frako) \leq 2\int_{V\times\calZ} \log \left(\frac{p_{\thetaMAP}(y,z)}{p_{\tilde{\theta}}(y,z)}\right)\mathrm{d}(\empM-\IP_{\theta_\frako})(y,z) + \frakd_\rate^2(\tilde{\theta},\theta_\frako)
	\end{equation}
	$\IP_{\theta_\frako}^N$-almost surely. Moreover, for any $R,\rate\in \pRZ$ and some fixed constant $\const_0\in\IR_{\geq 2}$ define the measurable event
	\begin{equation*}
		\Xi_{\rate,R} := \left\{\frakd_\rate^2(\thetaMAP,\theta_{\frako}) \geq \const_0\cdot(\frakd_\rate^2(\tilde{\theta},\theta_\frako)+R^2)\right\}.
	\end{equation*}
	On $\Xi_{\rate,R}$, we further have
	\begin{equation}\label{eq:MestimationII}
		\frakd_\rate^2(\thetaMAP,\tilde{\theta})\geq \const_1\cdot R^2
	\end{equation}
	for some $\const_1\in\IR_{\geq 1}$, as well as
	\begin{equation}\label{eq:MestimationIII}
		\frakd_\rate^2(\thetaMAP,\theta_\frako) -  \frakd_\rate^2(\tilde{\theta},\theta_\frako) \geq \frac{1}{6} \frakd_\rate^2(\thetaMAP,\tilde{\theta}).
	\end{equation}
\end{lemma}

\begin{proof}[Proof of \cref{lem: helpinglemmaI}]
	The proof of \cref{eq:MestimationI} is elementary and follows along the lines of Lemma 4.10 in \cite{NicklWang2023} (replacing the finite-dimensional approximation of $\theta_{\frako}$ by some arbitrary $\thetatilde\in\Theta$).
	We continue with the proof of \cref{eq:MestimationII}. On the event $\Xi_{\rate,R}$, we have
	\begin{align*}
		\const_0\cdot R^2 &\leq \frakd_\rate^2(\thetaMAP,\theta_\frako) - 	\const_0\cdot \frakd_\rate^2(\tilde{\theta},\theta_\frako)\\
		&\leq 2\left(\frac{1}{2}\norm{\scrG(\thetaMAP)-\scrG(\theta_\frako)}_{\IL^2_{\zeta}(\calZ,V)}^2-\norm{\scrG(\thetatilde)-\scrG(\thetatrue)}_{\IL^2_{\zeta}(\calZ,V)}^2\right)\\
		&\hspace{6cm}+\rate^2\left(\norm{\thetaMAP}_{H^\alpha(\domain,W)}^2-2\norm{\thetatilde}_{H^\alpha(\domain,W)}^2\right)\\
		&\leq 2\left(\norm{\scrG(\thetaMAP)-\scrG(\theta_\frako)}_{\IL^2_{\zeta}(\calZ,V)} -\norm{\scrG(\thetatilde)-\scrG(\thetatrue)}_{\IL^2_{\zeta}(\calZ,V)}\right)^2\\
		&\hspace{4cm}+ 2\rate^2\norm{\thetaMAP}_{H^\alpha(\domain,W)}^2 - \rate^2\norm{\thetaMAP}_{H^\alpha(\domain,W)}^2- 2\rate^2\norm{\thetatilde}_{H^\alpha(\domain,W)}^2\\
		&\leq 2\frakd_\rate^2(\thetaMAP,\thetatilde),
	\end{align*}
	where we have used $\frac{1}{2}a^2 - b^2 \leq (a-b)^2$ for two real numbers $a,b\in\IR$ as well as the reversed triangle inequality. Rearranging, the last inequality shows \cref{eq:MestimationII}. \cref{eq:MestimationIII} follows by the same elementary computations as in \cite{Geer_2001} (see pp. 3-4).
\end{proof}

	\begin{proof}[Proof of \cref{thm:maintheorem}]
		The proof follows the key ideas of \cite{BohrNickl2023}. Since we are allowing a general $\thetatilde\in\Theta$, we also need to adapt the ideas of \cite{NicklGeerWang2020}. We start with the exponential inequality \cref{eq:mainconcentration}. Fix $\thetatilde\in\Theta$. According to \cref{lem: helpinglemmaI}, it sufficies to control the $\IP_{\theta_\frako}^N$-probability of the event 
		\begin{equation*}
			\tilde{\Xi}_{\rate,R}:=\Xi_{\rate,R}\cap\left\{		\frakd_\rate^2(\thetaMAP,\theta_\frako) \leq 2\int_{V\times\calZ} \log \left(\frac{p_{\thetaMAP}(y,z)}{p_{\tilde{\theta}}(y,z)}\right)\mathrm{d}(\empM-\IP_{\theta_\frako})(y,z) + \frakd_\rate^2(\tilde{\theta},\theta_\frako)\right\}.
		\end{equation*}
		For any $\theta\in\Theta$, we define the empirical process
		\begin{equation*}
			\calW_N(\theta):= \int_{V\times\calZ} \log\left( \frac{p_{\theta}(y,z)}{p_{\tilde{\theta}}(y,z)}\right)\mathrm{d}(\empM-\IP_{\theta_\frako})(y,z) 
		\end{equation*}
		and obtain with \cref{eq:MestimationII} and \cref{eq:MestimationIII} from \cref{lem: helpinglemmaI}
		\magenta{
		\begin{align*}
			\IP_{\theta_\frako}^N(\tilde{\Xi}_{\rate,R}) 
			\leq\IP_{\theta_\frako}^N\left(|\calW_N(\thetaMAP)|\geq \frac{1}{12}\frakd_\rate^2(\thetaMAP,\thetatilde),\,\frakd_\rate^2(\thetaMAP,\thetatilde)\geq c_1^2 R^2,\right).
		\end{align*}}
		for some $\const_{1}\in\IR_{\geq1}$. In the following step, we apply a so-called \textit{peeling device}, where we refer to \cite{Geer2000} for a generic explanation. This method was already used in a similar way in \cite{NicklGeerWang2020} and \cite{BohrNickl2023} to control the process $\calW_N(\thetaMAP)$ over disjoint slices of $\Theta$. In detail, for all $l\in\IN$ and $\const_{1}$ sufficiently large, define the family of slices 
		\begin{equation*}
			\Theta_l:=\left\{\theta\in \Theta:\, 2^{2(l-1)}R^2\const_{1}^2 \leq \frakd_\rate^2(\theta,\thetatilde)< 2^{2l}R^2\const_{1}^2\right\} \subseteq H^{\alpha}(\domain, W,2^lR\rate^{-1}\const_{1}),
		\end{equation*}
		where the inclusion follows from the definition of $\frakd_\rate^2$. Note that $\{\theta\in\Theta:\,\frakd_\rate^2(\theta,\thetatilde) \geq R^2\const_1^2\} = \bigcup_{l\in\IN}\Theta_l$. We obtain
		\begin{equation*}
			\IP_{\theta_\frako}^N\left(\Xi_{\rate,R}\right) 
			\leq \sum_{l\in\IN}	\IP_{\theta_\frako}^N\left(\sup_{\theta\in\Theta_l}|\calW_{N}(\theta)|\geq \frac{2^{2l}}{48}R^2\const_{1}^2\right). 
		\end{equation*}
		We follow the approach of \cite{BohrNickl2023} and decompose $\calW_{N}$ into two empirical processes. Indeed, for each $\theta\in\Theta$, simple computations show 
		\begin{equation*}
			N \calW_{N} (\theta) 
			= \ell_N(\theta)-\ell_N(\thetatilde) - \IE_{\thetatrue}^N[\ell_N(\theta)-\ell_N(\thetatilde)],
		\end{equation*}
		which evidently defines a centered stochastic process. Given the data $D_N\sim\IP_{\theta_\frako}^N$, we further decompose the log-likelihood difference of the last display by
		\begin{align*}
			\ell_N(\theta)-\ell_N(\thetatilde) &=\frac{1}{2\sigma^2}\sum_{i\in\nset{N}}\left\{\Vnorm{Y_i-\scrG(\thetatilde)(Z_i)}^2-\Vnorm{Y_i-\scrG(\theta)(Z_i)}^2\right\}\\
			& \overset{d}{=}  \frac{1}{\sigma}\sum_{i\in\nset{N}}\langle \varepsilon_i,\scrG(\theta)(Z_i)-\scrG(\thetatilde)(Z_i) \rangle_V \\&\phantom{==}+  \frac{1}{2\sigma^2}\sum_{i\in\nset{N}}\left\{\Vnorm{\scrG(\theta_\frako)(Z_i)-\scrG(\thetatilde)(Z_i) }^2-\Vnorm{\scrG(\theta_\frako)(Z_i)-\scrG(\theta)(Z_i)}^2\right\}.
		\end{align*}
		Defining for all $\theta\in\Theta$
		\begin{equation}\label{eq:T1}
			T_1(\theta) :=\frac{1}{\sigma\sqrt{N}}\sum_{i\in\nset{N}} \langle \varepsilon_i,\scrG(\theta)(Z_i)-\scrG(\thetatilde)(Z_i) \rangle_V,
		\end{equation}
		and 
		\begin{equation}\label{eq:T3}
			T_2(\theta) :=\frac{1}{2\sigma^2\sqrt{N}}\sum_{i\in\nset{N}}\left\{\Vnorm{\scrG(\theta_\frako)(Z_i)-\scrG(\thetatilde)(Z_i) }^2-\Vnorm{\scrG(\theta_\frako)(Z_i)-\scrG(\theta)(Z_i)}^2\right\},
		\end{equation}
		respectively, we have that $	\IE_{\theta_\frako}^N[T_1(\theta)] =0$, as $Z$ and $\varepsilon$ are independent and $\varepsilon$ is centered. With these, we have derived so far
		\begin{align*}
			\IP_{\theta_\frako}^N\left(\Xi_{\rate,R}\right) &\leq \sum_{l\in\IN}	\IP_{\theta_\frako}^N\left(\sup_{\theta\in\Theta_l}|\calW_{N}(\theta)|\geq \frac{2^{2l}}{48}R^2\const_{1}^2\right)\\
			&\leq \sum_{l\in\IN}	\IP_{\theta_\frako}^N\left(\sup_{\theta\in\Theta_l}\left\{|T_1(\theta)|\right\}\geq \sqrt{N} \frac{2^{2l}}{96}R^2\const_{1}^2\right)\\ &\phantom{==}+\sum_{l\in\IN}\IP_{\theta_\frako}^N\left(\sup_{\theta\in\Theta_l}\left\{|T_2(\theta)-\IE_{\theta_\frako}^N[T_2(\theta)]|\right\}\geq \sqrt{N} \frac{2^{2l}}{96}R^2\const_{1}^2\right).
		\end{align*}
		We will discuss the two different probabilities separately and starting with the \textit{multiplier-type process} defined in \cref{eq:T1}. In order to apply \cref{lem:EPT}, we denote for each $j\in\nset{d_V}$ by $\scrG(\cdot)_j:\Theta\to\IL^2(\domain)$ the $j$-th coordinate of the forward map $\scrG:\Theta\to\IL^2(\domain,V)$, where we again identify $V$ with $\IR^{d_V}$. We define a class of measurable functions by 
		\begin{equation}\label{eq:functionclassH1}
			\calH_{l,j}:=\left\{\scrG(\thetatilde)_j-\scrG(\theta)_j:\calZ\to\IR:\, \theta\in\Theta_l\right\}.
		\end{equation}
		To shorten notation in the subsequent, we write $h_{\theta,j}:=\scrG(\thetatilde)_j-\scrG(\theta)_j$. Thus, for any $\theta\in\Theta$, we have
		\begin{align*}
			T_1(\theta)  %&= - \frac{1}{\sqrt{N}}\sum_{i\in\nset{N}}\Vpro{\scrG(\theta_\frako)(Z_i)-\scrG(\theta)(Z_i)}{\varepsilon_i} \\
			& \overset{d}{=} - \frac{1}{\sigma\sqrt{N}} \sum_{i\in\nset{N}}\sum_{j\in\nset{d_V}}(\scrG(\thetatilde)_j(Z_i)-\scrG(\theta)_j(Z_i))\varepsilon_{ij}\\
			&= -  \frac{1}{\sigma\sqrt{N}}\sum_{i\in\nset{N}}\sum_{j\in\nset{d_V}}h_{\theta,j}(Z_i)\varepsilon_{ij}
			=: - \frac{1}{\sigma}\sum_{j\in\nset{d_V}} Z_{N,j}^{(1)}(\theta),
		\end{align*}
		where $\varepsilon_{ij}\overset{i.i.d.}{\sim}\operatorname{N}(0,1)$ for $(i,j)\in\nset{N}\times\nset{d_V}$. Conclusively, we see
		\begin{equation*}
			\IP_{\theta_\frako}^N\left(\sup_{\theta\in\Theta_l}|T_1(\theta)|\geq \sqrt{N} \frac{2^{2l}}{96}R^2\const_{1}^2\right)\leq d_V\max_{j\in\nset{d_V}}\left\{	\IP_{\theta_\frako}^N\left(\sup_{\theta\in\Theta_l}|Z_{N,j}^{(1)}(\theta)|\geq \sqrt{N} \frac{2^{2l}}{96d_V}R^2\const_{1}^2\sigma\right)\right\}.
		\end{equation*}
		We now verify that the function class $\calH_{l,j}$ defined in \cref{eq:functionclassH1} satisfies the assumptions formulated in \cref{lem:EPT}. Indeed, by C2) in \cref{ass:ForwardMapI}, we have uniformly for all $\theta\in\Theta_l$
		\begin{equation*}
			\sup_{z\in\calZ}\left\{|h_{\theta,j}(z)|\right\} = \sup_{z\in\calZ}\left\{|\scrG(\thetatilde)_j(z)-\scrG(\theta)_j(z)|\right\} \leq  \sup_{z\in\calZ}\left\{\Vnorm{\scrG(\thetatilde)(z)-\scrG(\theta)(z)}\right\} \leq 2\Const_{\scrG,\operatorname{B}}.
		\end{equation*}
		Moreover, for $Z\sim\zeta$ and all $\theta\in\Theta_l$, $l\in\IN$, we have
		\begin{equation*}
			\IE^Z[h_{\theta,j}^2(Z)] %= \int_\calZ \Vnorm{\scrG(\theta)_j(Z)-\scrG(\theta_\frako)_j(Z)}^2\mathrm{d}\zeta(z)
			\leq \sum_{j\in\nset{d_V}} \int_\calZ \Vnorm{\scrG(\theta)_j(z)-\scrG(\thetatilde)_j(z)}^2\mathrm{d}\zeta(z)
			= \norm{\scrG(\theta)-\scrG(\thetatilde)}_{\IL^2_{\zeta}(\calZ,V)}^2.
		\end{equation*}
		Observing that $\theta\in\Theta_l$, the last expression is further upper bounded by 
		\begin{equation}\label{eq:v_l}
			\norm{\scrG(\theta)-\scrG(\thetatilde)}_{\IL^2_{\zeta}(\calZ,V)}^2\leq \frakd_\rate^2(\theta,\thetatilde) \leq 2^{2l}R^2 \const_{1}^2=:  v_l^2.
		\end{equation}
		Additionally, we need to upper bound the $\IL^2$- and $\IL^\infty$-metric entropy integrals $\calJ_2$ and $\calJ_\infty$. Thus, we proceed by upper bounding the respective metric entropies
		\begin{equation*}
			\log N\left(\calH_{l,j},d_2,\rho\right)	\text{ and }\log N\left(\calH_{l,j},d_\infty,\rho\right),
		\end{equation*}
		where $d_2$ and $d_\infty$ are defined according to \cref{lem:EPT}. Starting with the first one, let $\theta_1,\theta_2\in\Theta_l$ be arbitrarily chosen. Exploiting C1) in \cref{ass:ForwardMapI}, we have 
		\begin{equation*}
			d_2(\theta_1,\theta_2) = \norm{h_{\theta_1,j}-h_{\theta_2,j}}_{\IL^2_\zeta(\calZ)} \hspace{-0.5ex}
			 \leq\Const_{\Lip,2}\left(1+\norm{\theta_1}^{\gamma_1}_{H^\alpha(\domain,W)}\hspace{-1ex}\lor\norm{\theta_2}^{\gamma_1}_{H^\alpha(\domain,W)}\right)\hspace{-0.5ex}\norm{\theta_1-\theta_2}_{(H^\kappa(\domain,W))^*}.
		\end{equation*}
		By definition of $\Theta_l$, $\theta_1$, $\theta_2$ are contained in a $H^\alpha$-ball, such that the polynomial growth in the last display can be upper bounded by the corresponding radius. Hence, it suffices to find a $\rho$-covering of $\Theta_l$ with respect to the dual $(H^\kappa(\domain,W))^*$-norm. Applying \cref{lem:SobolevCovering} i), there exists a constant $\const_{2}=\const_{2}(\alpha,\kappa,d)\in\pRZ$, such that for all $\rho\in\pRZ$
		\begin{equation*}
			N(\Theta_l,\norm{\cdot}_{(H^\kappa(\domain,W))^*},\rho) 
			\leq \exp\left(\const_{2} \cdot d_W \left( \frac{2^lR\sqrt{d_W}\const_1}{\rate\rho}\right)^{\frac{d}{\alpha+\kappa}}\right).
		\end{equation*}
		 As a conseqence, for every $\rho\in\pRZ$, there exists a $\rho$-covering $\{\theta^{(i)}\}_{i\in\nset{N}}$ of $\Theta_l$ \magenta{with respect to $\norm{\cdot}_{(H^\kappa(\domain,W))^*}$} with $N=N(\Theta_l,\norm{\cdot}_{(H^\kappa(\domain,W))^*},\rho)$. Hence, for all $\theta\in\Theta_l$, there exists some $i_\circ\in\nset{N}$ with 
		\begin{align*}
			d_2(\theta,\theta^{(i_\circ)}) = \norm{h_{\theta,j}-h_{\theta^{(i_\circ)},j}}_{\IL^2_\zeta(\calZ)}
			& \leq  \Const_{\Lip,2}\left(1+\left(\frac{2^lR\const_{1}}{\rate}\right)^{\gamma_1}\right)\cdot\rho
			=: m_1\cdot \rho
		\end{align*}
		and we deduce
		\begin{equation*}
			\log N\left(\calH_{l,j},d_2,\rho\right) \leq \const_{2} \cdot d_W \left( \frac{2^lR\sqrt{d_W}m_1\const_1}{\rate\rho}\right)^{\frac{d}{\alpha+\kappa}} = \const_{3} \cdot  \left( \frac{2^lRm_1\const_1}{\rate\rho}\right)^{\frac{d}{\alpha+\kappa}} 
		\end{equation*}
		for some constant $\const_{3} = \const_{3}(\alpha,\kappa,d,d_W)\in\pRZ$. An analogous argument yields a similar upper bound for the $d_\infty$-metric entropy. Indeed, for any $\theta_1,\theta_2\in\Theta_l$, exploiting C3) in \cref{ass:ForwardMapI}, we have 
		\begin{equation*}
			d_\infty(\theta_1,\theta_2) = \norm{h_{\theta_1,j}-h_{\theta_2,j}}_\infty\leq\Const_{\operatorname{Lip},\infty}\left(1+\norm{\theta_1}_{C^{\eta_1}(\domain,W)}^{\gamma_2}\lor\norm{\theta_2}_{C^{\eta_1}(\domain,W)}^{\gamma_2}\right) \norm{\theta_1-\theta_2}_{C^{\eta_2}(\domain,W)}.
		\end{equation*}
		Using that $\alpha >\eta_1\lor\eta_2+\frac{d}{2}$, the polynomial growth can be bounded due to the continuous embedding $H^\alpha(\domain,W)\hookrightarrow C^{\eta_1}(\domain,W)$. Hence, it suffices to find a covering for $\Theta_l$ with respect to the $C^{\eta_2}$-norm. By an application of \cref{lem:SobolevCovering} ii), we have 
		\begin{equation*}
		N(\Theta_l,\norm{\cdot}_{C^{\eta_2(\domain,W)}},\rho) \leq \exp\left(\const(\alpha,\eta_2,d) \cdot d_W \left( \frac{2^lR\sqrt{d_W}\const_1}{\rate\rho}\right)^{\frac{d}{\alpha-\eta_2}}\right).
		\end{equation*}
		 \magenta{As a conseqence, for every $\rho\in\pRZ$, there exists a $\rho$-covering $\{\theta^{(i)}\}_{i\in\nset{N}}$ of $\Theta_l$ with respect to $\norm{\cdot}_{C^{\eta_2(\domain,W)}}$ with $N=N(\Theta_l,\norm{\cdot}_{C^{\eta_2(\domain,W)}},\rho)$.} Thus, for all $\theta\in\Theta_l$, there exists some $i_\circ\in\nset{N}$ with
			\begin{align*}
				d_\infty(\theta,\theta^{(i_\circ)}) \leq \magenta{\norm{h_{\theta,j}-h_{\theta^{(i_\circ)},j}}_\infty}\leq \const(\Const_{\Lip,\infty},\gamma_2)\cdot\left(1+\left(\frac{2^l R \const_1}{\rate}\right)^{\gamma_2}\right)\cdot \rho\magenta{=:m_2\cdot\rho}
			\end{align*}
			\magenta{and we deduce}
			\begin{equation*}
				\log N\left(\calH_{l,j},d_\infty,\rho\right) \leq \const_{4}\cdot\left(\frac{2^lRm_2\const_1}{\rate\rho}\right)^{\frac{d}{\alpha-\eta}}
			\end{equation*}
			\magenta{for some constant $\const_{4} = \const_{4}(\alpha,\eta_2,d,d_W)\in\pRZ$.} We are now able to upper bound the corresponding entropy integrals. Indeed, as $2(\alpha+\kappa)>d$, we have
		\begin{align*}
			\calJ_2(\calH_{l,j})  &= \int_0^{4 v_l}\sqrt{\log N\left(\calH_{l,j},d_2,\rho\right)}\mathrm{d}\rho\\
			&\leq \sqrt{\const_{3}}\cdot\int_0^{4 v_l}\left( \frac{2^lRm_1\const_1}{\rate\rho}\right)^{\frac{d}{2(\alpha+\kappa)}} \mathrm{d}\rho
			%& \leq \int_0^{4\sigma_s} \const_{d_W,d,\alpha,\kappa}\left(\frac{2^lRm_1}{\rate\rho}\right)^{\frac{d}{2(\alpha+\kappa)}} \mathrm{d}\rho \\
			%& =  \const_{d_W,d,\alpha,\kappa}\left(\frac{2^lRm_1}{\rate}\right)^{\frac{d}{2(\alpha+\kappa)}}\int_0^{4\sigma_s}\left(\frac{1}{\rho}\right)^{\frac{d}{2(\alpha+\kappa)}} \mathrm{d}\rho \\
			=  \sqrt{\const_{3}}\left(\frac{2^lRm_1\const_1}{\rate}\right)^{\frac{d}{2(\alpha+\kappa)}}  v_l^{1-{\frac{d}{2(\alpha+\kappa)}}},
		\end{align*}
		where the numerical value of constant $\const_3$ changed without any further dependencies. To shorten notation in the further analysis write $\tau_1 :=\frac{d}{2(\alpha+\kappa)}\in(0,1)$. By assumption $R\geq\rate$. We have
		\begin{align*}
			\left(\frac{2^lRm_1\const_1}{\rate}\right)^{\frac{d}{2(\alpha+\kappa)}}  v_l^{1-{\frac{d}{2(\alpha+\kappa)}}}\hspace{-1ex} & = c_1^{\tau_1}\cdot \Const_{\Lip,2}^{\tau_1}\cdot v_l^2\cdot\left(1+\left(\frac{2^lR\const_{1}}{\rate}\right)^{\gamma_1}\right)^{\tau_1}\cdot2^{l\tau_1}\cdot R^{\tau_1}\cdot\rate^{-\tau_1}\cdot  v_l^{-1-\tau_1}\\
			&  \leq  c_1^{\tau_1+\gamma_1\tau_1}\cdot \Const_{\Lip,2}^{\tau_1}\cdot v_l^2\cdot 2^{\tau_1}\cdot 2^{l\tau_1+l\tau_1\gamma_1}\cdot R^{\tau_1+\gamma_1\tau_1}\cdot\rate^{-\tau_1-\tau_1\gamma_1}\cdot  v_l^{-1-\tau_1}\\
			& = c_1^{-(1-\gamma_1\tau_1)}\cdot \Const_{\Lip,2}^{\tau_1}\cdot v_l^2\cdot 2^{\tau_1}\cdot 2^{-l(1-\tau_1\gamma_1)}\cdot R^{-1+\gamma_1\tau_1}\cdot\rate^{-\tau_1(+\tau_1\gamma_1)}\\
			&\leq c_1^{-(1-\gamma_1\tau_1)}\cdot \Const_{\Lip,2}^{\tau_1}\cdot v_l^2\cdot 2^{\tau_1}\cdot \sqrt{N},
		\end{align*}
		where we have used \cref{eq:main:ass1} and the fact that $\alpha+\kappa>\frac{d}{2}\gamma_1$. In particular, the last bound implies the existence of a constant $c_5 =c_5(\alpha,\gamma_1,\kappa,d,d_W, \Const_{\Lip,2})\in\pRZ$, such that
		\begin{equation*}
			\calJ_2(\calH_{l,j}) \leq c_5\cdot c_1^{-(1-\gamma_1\tau_1)}\cdot v_l^2\cdot  \sqrt{N}.
		\end{equation*}
		Analogously, as $\alpha-\eta_2 > d$, we have
		\begin{align*}
			\calJ_\infty(\calH_{l,j}) &= \int_0^{8\Const_{\scrG,\operatorname{B}}}\log N\left(\calH_{l,j},d_\infty,\rho\right)\mathrm{d}\rho \\
			&\leq\const_{4}\cdot\int_0^{8\Const_{\scrG,\operatorname{B}}}\left(\frac{2^lRm_2c_1}{\rate\rho}\right)^{\frac{d}{\alpha-\eta_2}}\mathrm{d}\rho
			\leq \const(\alpha,\eta_2,d,d_W,\Const_{\scrG,\operatorname{B}})\cdot\left(\frac{2^lRm_2c_1}{\rate}\right)^{\frac{d}{\alpha-\eta_2}}.
		\end{align*}
		Using \cref{eq:main:ass2} with similar arguments as above, one can verify analogously
		\begin{equation*}
			\calJ_\infty(\calH_{l,j}) \leq\const_{6}  v_l^2 N
		\end{equation*}
		with some $\const_{6} =\const_{6}(\alpha,\eta_1,\eta_2,\gamma_2,\Const_{\Lip,\infty},\Const_{\scrG,\operatorname{B}},d,d_W,\const_1)\in\pRZ$, which can be also as small as desired by choosing $c_1$ sufficiently large. Now in the setting of \cref{lem:EPT}, for arbitrary $\bar{\const}\in\IR_{\geq1}$, such that $x:=\bar{\const}2^{2l}NR^2\in\IR_{\geq 1}$, we have for some \magenta{universal} constant $L\in\pRZ$
		\magenta{
			\begin{align*}
				&L\left(\calJ_2(\calH_{l,j})+ v_l\sqrt{x}+N^{-\frac{1}{2}}(\calJ_\infty(\calH_{l,j})+8\Const_{\scrG,\operatorname{B}}x)\right)\\
				&\hspace{1cm}\leq L\left(c_5c_1^{-(1-\gamma_1\tau_1)}v_l^2\sqrt{N}+v_l\sqrt{\bar{c}}\sqrt{N}2^lR+\frac{1}{\sqrt{N}}(c_6v_l^2N+8\bar{c}N2^{2l}R^2\Const_{\scrG,\operatorname{B}})\right)\\
				&\hspace{1cm}\leq \sqrt{N} \frac{ v_l^2\sigma}{96\cdot d_V}\left(96Ld_V\sigma^{-1}c_5c_1^{-(1-\gamma_1\tau_1)}+96Ld_V\frac{\sqrt{\bar{c}}}{\sigma c_1}+96Ld_Vc_6+768Ld_V\Const_{\scrG,\operatorname{B}}\frac{\bar{c}}{c_1^2}\right)\\
				&\hspace{1cm}\leq \sqrt{N} \frac{ v_l^2\sigma}{96\cdot d_V},
			\end{align*}
			where in the last step follows after choosing  $\const_{1}=\const_{1}(\texttt{P})\in\pRZ$ sufficiently large.} Hence an application of \cref{lem:EPT} yields
		\begin{align*}
			&\IP_{\theta_\frako}^N\left(\sup_{\theta\in\Theta_l}|Z_{N,j}^{(1)}(\theta)|\geq \sqrt{N} \frac{ v_l^2\sigma}{96d_V}\right)\\
			&
			\hspace{1.5cm}\leq \IP_{\theta_\frako}^N\left(	\sup_{\theta\in\Theta_l}|Z_{N,j}^{(1)}(\theta)|\geq L\left(\calJ_2(\calH_{l,j})+ v_l\sqrt{x}+N^{-\frac{1}{2}}(\calJ_\infty(\calH_{l,j})+8\Const_{\scrG,\operatorname{B}}x)\right) \right)\\
			&\hspace{1.5cm}\leq 2\exp\left(-\bar{\const}  2^{2l}NR^2\right).
		\end{align*}
		As the last upper bound is independent of the choice of $j\in\nset{d_V}$, we have overall
		\begin{equation*}
			\IP_{\theta_\frako}^N\left(\sup_{\theta\in\Theta_l}\left\{|T_1(\theta)|\right\}\geq \sqrt{N} \frac{ v_l^2}{96}\right)\leq 2d_V\exp\left(-\bar{\const} 2^{2l}NR^2\right).
		\end{equation*}
		It remains to study the process $T_2$. By writing 
		\begin{equation*}
				T_2(\theta)
				=\frac{1}{2\sigma^2\sqrt{N}}\sum_{i\in\nset{N}}h_{\theta}(Z_i)
		\end{equation*}
		with $h_{\theta}(\cdot) = \Vnorm{\scrG(\theta_\frako)(\cdot)-\scrG(\thetatilde)(\cdot) }^2-\Vnorm{\scrG(\theta_\frako)(\cdot)-\scrG(\theta)(\cdot)}^2$, we need to control
		\begin{equation*}
			\IP_{\theta_\frako}^N\left(\sup_{\theta\in\Theta_l}\left\{\left|\frac{1}{\sqrt{N}}\sum_{i\in\nset{N}}(h_{\theta}(Z_i)-\IE_{\theta_\frako}^N[h_{\theta}(Z_i)])\right|\right\}\geq \sqrt{N} \frac{2^{2l}}{48}R^2\const_{1}^2\sigma^2\right).
		\end{equation*}
		Similiar as before, we define for each $l\in\IN$ a class of measurable functions
		\begin{equation*}
			\calH_l :=\left\{h_\theta:=\Vnorm{\scrG(\theta_\frako)(\cdot)-\scrG(\thetatilde)(\cdot) }^2-\Vnorm{\scrG(\theta_\frako)(\cdot)-\scrG(\theta)(\cdot)}^2:\calZ\to\IR\,:\,\theta\in\Theta_l\right\},
		\end{equation*}
		and verify the assumptions of \cref{lem:EPT}. Indeed, note that for each $\theta\in\Theta_l$, we have 
		\begin{align*}
			\sup_{z\in\calZ}\left\{\left|h_\theta(z)\right|\right\} &= 	\sup_{z\in\calZ}\left\{\left|\Vnorm{\scrG(\theta_\frako)(z)-\scrG(\thetatilde)(z) }^2-\Vnorm{\scrG(\theta_\frako)(z)-\scrG(\theta)(z)}^2\right|\right\} \leq 8\Const_{\scrG,\operatorname{B}} 
		\end{align*}
		due to C2) in \cref{ass:ForwardMapI}. Moreover, for $Z\sim\zeta$ and all $\theta\in\Theta_l$, $l\in\IN$, we have 
		\begin{align*}
			\IE^Z[h_{\theta}^2(Z)] & = \IE^Z\left[\left|\Vnorm{\scrG(\theta_\frako)(Z)-\scrG(\thetatilde)(Z) }^2-\Vnorm{\scrG(\theta_\frako)(Z)-\scrG(\theta)(Z)}^2\right|^2\right]\\
			& \leq 16\Const_{\scrG,\operatorname{B}}^2\norm{\scrG(\thetatilde)-\scrG(\theta)}_{\IL^2_\zeta(\calZ,V)}^2 \leq 16\Const_{\scrG,\operatorname{B}}^2v_l^2,
		\end{align*}
		with $v_l$ defined in \cref{eq:v_l}, where we have applied C2) from \cref{ass:ForwardMapI}. In order to bound the entropy integrals, observe that for any $\theta_1,\theta_2\in\Theta_l$ 
		\begin{align*}
			d_2(\theta_1,\theta_2)  &= \norm{h_{\theta_1}-h_{\theta_2}}_{\IL^2_\zeta(\calZ)}
			 = \norm{\Vnorm{\scrG(\theta_\frako)(\cdot)-\scrG(\theta_2)(\cdot)}^2-\Vnorm{\scrG(\theta_\frako)(\cdot)-\scrG(\theta_1)(\cdot)}^2}_{\IL^2_\zeta(\calZ)}\\
			&\leq 4\Const_{\scrG,\operatorname{B}}\cdot\norm{\scrG(\theta_1)-\scrG(\theta_2)}_{\IL^2_\zeta(\calZ,V)},
		\end{align*}
		as well as
		\begin{align*}
			d_\infty(\theta_1,\theta_2)  &= \norm{h_{\theta_1}-h_{\theta_2}}_{\infty}
			 = \norm{\Vnorm{\scrG(\theta_\frako)(\cdot)-\scrG(\theta_2)(\cdot)}^2-\Vnorm{\scrG(\theta_\frako)(\cdot)-\scrG(\theta_1)(\cdot)}^2}_{\infty}\\
			&\leq 4\Const_{\scrG,\operatorname{B}}\cdot\norm{\scrG(\theta_1)-\scrG(\theta_2)}_{\infty}.
		\end{align*} 
		Evidently, we obtain the same upper bounds for the metric entropies
	\begin{equation*}
		\log N\left(\calH_{l},d_2,\rho\right)	\text{ and }\log N\left(\calH_{l},d_\infty,\rho\right)
	\end{equation*}
		 as already derived before (up to a change of constants). Thus, following the same argumentation as above, we can choose $c_1 =c_1(\texttt{P},d_V)\in\pRZ$ sufficiently large to obtain
		\begin{equation*}
			\IP_{\theta_\frako}^N\left(\sup_{\theta\in\Theta_l}\left\{|T_2(\theta)-\IE_{\theta_\frako}^N[T_2(\theta)]|\right\}\geq \sqrt{N} \frac{2^{2l}}{96}R^2\const_{1}^2\right) \leq 2\exp\left(-\bar{\const}2^{2l}NR^2\right).
		\end{equation*}	
		Overall, we have derived the existence of constants $\const_{7} = \const_{7}(d_V)\in\pRZ$, such that
		\begin{equation*}
			\IP_{\theta_\frako}^N(\Xi_{\rate,R}) \leq \const_{7}\cdot\sum_{l\in\IN}\exp(-\bar{\const}2^{2l}NR^2) \leq \const_{7}\exp(-\bar{\const}NR^2),
		\end{equation*}
		if $c_1\in\IR_{\geq 1}$ in the definition of $\Xi_{\rate,R}$ is choosen sufficiently large, observing that the numerical values of the constants may changed again in the last line without any further dependencies. 
		This shows the exponential inequality in i). Now we verify the upper bound of the expectation, i.e. ii). Let $\thetatilde\in\Theta$ be fixed. For some appropriate $\frakm\in\pRZ$ specified later, we have
		\begin{align*}
			\IE_{\theta_\frako}^N[\frakd_\rate^2(\thetaMAP,\theta_{\frako})] & \leq \frakm +\int_{(\frakm,\infty)} 	\IP_{\theta_\frako}^N(\frakd_\rate^2(\thetaMAP,\theta_{\frako})\geq t)\mathrm{d}\calL(t).
		\end{align*}
		Let $\rate\in\pRZ$ satisfy \cref{eq:main:ass3} and fix some $\bar{\const}\in\IR_{\geq 1}$. We define the map $	T:(\rate^2,\infty)\to\pRZ,\,	t\mapsto \const_{1}\cdot(\frakd_\rate^2(\thetatilde,\theta_{\frako})+t)$ and set $\frakm := \const_1\cdot(\frakd_\rate^2(\thetatilde,\theta_{\frako})+\rate^2)$, with $\const_1\in\IR_{>0}$ choosen sufficiently large according to i). Applying the transformation theorem, we obtain
		\begin{align*}
			\IE_{\theta_\frako}^N[\frakd_\rate^2(\thetaMAP,\theta_{\frako})]  & \leq \frakm +\const_1\cdot\int_{(\rate^2,\infty)} \IP_{\theta_\frako}^N\left(\frakd_\rate^2(\thetaMAP,\theta_{\frako})\geq \const_1\cdot(\frakd_\rate^2(\thetatilde,\theta_{\frako})+t)\right)\mathrm{d}\calL(t).
		\end{align*}
		Now observe that $t\in(\rate^2,\infty)$ satisfies $\sqrt{t}\geq \rate$ and hence also \cref{eq:main:ass1} and \cref{eq:main:ass2} by an application of \cref{eq:main:ass3}. Thus, utilizing the concentration inequality yields
		\begin{align*}
			\IE_{\theta_\frako}^N[\frakd_\rate^2(\thetaMAP,\theta_{\frako})]  & \leq \frakm +\const_1\cdot c(d_V)\cdot\int_{(\rate^2,\infty)} \exp(-\bar{\const}Nt)\mathrm{d}\calL(t) \leq \frakm+\frac{\const_{1}\cdot c(d_V)}{\bar{\const}N}.
		\end{align*}
		 The claim then follows by taking the infimum over all admissible $\thetatilde$.
	\end{proof}
\begin{proof}[Proof of \cref{co:MAPConsistency}]
	To prove the claim, we aim to apply \cref{thm:maintheorem}. Therefore choose $\thetatilde = \theta_{\frako}\in\Theta\cap\scrH(\operatorname{L})$. Let $\rate =\rate_{N}$ be as in the hypothesis. Elementary computations show that \cref{eq:main:ass1}, \cref{eq:main:ass2} and \cref{eq:main:ass3} are verified for the choice $R=\rate_{N}$. Moreover, observe for any $\const_1\in\pRZ$ that
	\begin{equation*}
			\const_1\cdot(\dist_\rate^2(\theta_{\frako},\theta_\frako)+R^2) = \const_1\cdot(\rate^2\norm{\theta_{\frako}}_{H^\alpha(\domain,W)}^2+\rate^2) \leq \const_1\cdot(\operatorname{L}^2+1)\cdot \rate^2.
	\end{equation*}
	Thus, for any $\bar{\const}\in\IR_{\geq 1}$, we can choose $\const_{1}=c_1(\texttt{P},L)\in\IR_{\geq 1}$ sufficiently large, that 
	\begin{multline}
		\IP_{\theta_\frako}^N\left(\dist_\rate^2(\thetaMAP,\theta_{\frako}) \geq  \const_1\cdot(\operatorname{L}^2+1)\cdot \rate^2\right)\\\leq 	\IP_{\theta_\frako}^N\left(\dist_\rate^2(\thetaMAP,\theta_{\frako}) \geq 2(\dist_\rate^2(\tilde{\theta},\theta_\frako)+R^2)\right)
		\leq \const(d_V)\cdot \exp(-\bar{\const} N\rate_N^2).
	\end{multline}
	This shows the first part of the claim. Moreover, according to the second part of \cref{thm:maintheorem}, we have
	\begin{equation*}
		\IE_{\theta_0}^N\left[\norm{\scrG(\thetaMAP)-\scrG(\theta_\frako)}^2_{\IL^2_\zeta(\calZ,V)}\right] \leq \const_{3}\left( \dist_\rate^2(\theta_{\frako},\theta_{\frako}) + \rate_N^2 + \frac{1}{N}\right)
		\leq\const_{3}\left( (\operatorname{L}^2+1)\rate_N^2 +  \frac{1}{N}\right)
		\leq \const_4\rate_N^2
	\end{equation*}
	for some constant $\const_{4}=\const_{4}(\operatorname{L},\alpha,\eta_1,\eta_2,\kappa,\gamma_1,\gamma_2,\Const_{\Lip,2},\Const_{\Lip,\infty},\Const_{\scrG,\operatorname{B}},d,d_W)\in\pRZ$. Taking the supremum over all admissible $\theta_{\frako}$ on both sides shows the claim.
\end{proof}
\begin{proof}[Proof of \cref{cor:stability}]
	\magenta{
		Firstly, for $c_1\in\IR_{>0}$ define the event 
		\begin{equation*}
			A_N :=\left\{\norm{\scrG(\thetaMAP)-\scrG(\theta_\frako)}_{\IL^2_\zeta(\calZ,V)}^2+\rate_N^2\norm{\thetaMAP}^2_{H^\alpha(\domain,W)} \leq  \const_{1} \rate_N^2\right\}.
		\end{equation*}
		By an application of \cref{co:MAPConsistency}, we can choose for each $\bar{\const}\in\IR_{\geq 1}$ the constant $\const_1=\const_1(\bar{\const},,\alpha,\eta_1,\eta_2,\kappa,\gamma_1,\gamma_2,\Const_{\Lip,2},\Const_{\Lip,\infty},\Const_{\scrG,\operatorname{B}},d,d_W,\norm{\theta_{\frako}}_{H^\alpha(\domain,W)})\in\IR_{\geq 1}$ large enough to obatin
		\begin{equation*}
			\IP_{\theta_\frako}^N\left(A_N\right) \geq 1- \const(d_V)\cdot\exp(-\bar{\const}N\rate_N^2).
		\end{equation*}
		We observe that 
		\begin{equation*}
			A_N \subset B_N:= \left\{ \norm{\scrG(\thetaMAP)-\scrG(\theta_\frako)}_{\IL^2_\zeta(\calZ,V)}^2 \leq \const_{1} \rate_N^2,\,\norm{\thetaMAP}_{H^\alpha(\domain,W)} \leq \sqrt{\const_{1}}\right\}.
		\end{equation*}
		As $\theta_{\frako}\in H^\alpha(\domain,W)$ is fixed, there exists some constant $\operatorname{L}=\operatorname{L}(\theta_\frako)\in\pRZ$ with $\norm{\theta_\frako}_{H^\alpha(\domain,W)}\leq\operatorname{L}$ and thus on the event $B_N$, we have
		\begin{equation*}
			\norm{\scrG(\thetaMAP)-\scrG(\theta_\frako)}_{\IL^2_\zeta(\calZ,V)}^2 \leq  \const_1\cdot\rate_N^2\text{ as well as } 	\norm{\thetaMAP}_{H^\alpha(\domain,W)}  +\norm{\theta_\frako}_{H^\alpha(\domain,W)} \leq \sqrt{\const_1}+L.
		\end{equation*}
		By applying the inverse continuity modulus formulated in \cref{ass:ForwardMapII}, we have on the event $B_N$
		\begin{equation*}
			\norm{\thetaMAP-\theta_\frako}_{\IL^2(\domain,W)}\leq \const(\Const_{\operatorname{M}},\const_1)\cdot \rate_N^\tau.
		\end{equation*}
		With the estimations above, the claim follows. 
	}
\end{proof}

%% file: content/appendix/app_pde.tex
\subsection{Proofs of \cref{sec:PDE}}\label{subsec:app:PDE}
In the following, we prove the results of \cref{sec:PDE}. The goal is to verify that the forward map $\scrG^D$ of the Darcy problem satisfies the conditions formulated in \cref{ass:ForwardMapI} and to apply the results of \cref{sec:model} afterwards.
 In the subsequent, we denote by $\Zygmund{s}{\domain} = B_{\infty\infty}^s(\domain)$ the Hölder-Zygmund space for which we have $\Zygmund{s'}{\domain} \hookrightarrow C^s(\domain)\hookrightarrow \Zygmund{s}{\domain}$, if $s'>s>0$ (see \cite{Triebel1983} \magenta{Section 3.4.2} for definitions and properties). 
Let $f\in C^1(\bar{\domain})$ with $f\geq f_{\min}\in\pRZ$. We define the differential operator 
\begin{equation*}
	L^D_f:	H_0^2(\domain) \to \IL^2(\domain)\text{, }	u \mapsto L_f^D[u] := \divergence(f\nabla u),
\end{equation*}
where $\nabla$ denotes the usual gradient and $\divergence$ the divergence operator. \magenta{From classical elliptic PDE theory (see for instance \cite{Evans2010} Theorem 4 and the Remark after it in Chapter 6.3), it is well-known, that $L^D$ is a linear and bounded operator, which possess a linear inverse, we denote by}
\begin{equation*}
	V_f^D:\IL^2(\domain)\to H^2_0(\domain)\text{, }h\mapsto V_f^D[h].
\end{equation*}
Thus, for any $h\in\IL^2(\domain)$, there exists a unique weak solution $w_{f,h}:=V^D_f[h]\in H^2_0(\domain)$ solving the Dirichlet boundary value problem 
\begin{equation*}
	\begin{dcases}
		L_f[w_{f,h}] &= h \text{ on $\domain$,}\\
		\hfill w_{f,h} &= 0  \text{ on $\partial\domain$.}
	\end{dcases}
\end{equation*} 
Utilizing the corresponding Feynman-Kac formula, Lemma 20 in \cite{NicklGeerWang2020} demonstrates that
\begin{equation}\label{eq:darcy:inftynormbound}
	\norm{V^D_f[h]}_\infty \leq \const(d,\domain,f_{\min})\cdot\norm{h}_\infty
\end{equation}
\magenta{holds} for all $h\in C^\eta(\domain)$, $\eta\in\pRZ$, if additionally $f\in C^2(\domain)$. Keeping the notation of \cref{sec:PDE} in mind, we now verify that the forward map $\scrG^D$ satisfies the conditions C1) - C3) formulated in \cref{ass:ForwardMapI}. \magenta{The cruicial part is to verify condition C3) for some $\eta_1$, $\eta_2$ and $\gamma_2$.}
\begin{enumerate}
	\item[C1)] $\IL^2$-Lipschitz continuity: Let $\alpha > 2 + \frac{d}{2}$. The forward map $\scrG^D:H^\alpha_c(\domain)\to\IL^2(\domain)$ is $(\alpha,4,1)$-regular, that is for all $\theta_1,\theta_2\in H^\alpha_c(\domain)$, we have
	\begin{align*}
		\norm{\scrG^D(\theta_1)-\scrG^D(\theta_2)}_{\IL^2(\domain)}\leq \Const_{\Lip,2}\cdot(1+\norm{\theta_1}_{H^\alpha(\domain)}^4\lor \norm{\theta_2}_{H^\alpha(\domain)}^4)\cdot\norm{\theta_1-\theta_2}_{(H^1(\domain))^*}
	\end{align*}
	with some universal constant $\Const_{\Lip,2}\in\pRZ$. The proof can be found in \cite{NicklGeerWang2020}.
	\item[C2)] Uniform boundedness: As described in \cref{sec:PDE}, we identify with each $\theta\in H^\alpha_c(\domain)$ some $f=f_\theta=\Psi(\theta)\in\calF$, such that $\scrG^D(\theta) = G^D(f_\theta) = V_{f_\theta}^D[g]$. Thus, with \cref{eq:darcy:inftynormbound}, we obtain
	\begin{align*}
		\sup_{\theta\in H^\alpha_c(\domain),x\in\domain}\left\{|\scrG^D(\theta)(x)|\right\} &=	\sup_{f\in \calF,x\in\domain}\left\{|G^D(f)(x)|\right\} \\
		 &=	\sup_{f\in \calF,x\in\domain}\left\{|V_f^D[g](x)|\right\} 
		 \leq \const(d,\domain,f_{\min})\cdot \norm{g}_\infty.
	\end{align*}
	
	\item[C3)] $\IL^\infty$-Lipschitz continuity: Let $\alpha>2+\frac{d}{2}$. Let $\theta_1,\theta_2\in H^\alpha_c(\domain)$ and identify $f_{\theta_1}=\Psi(\theta_1),f_{\theta_2}=\Psi(\theta_2)\in\calF$. With above notation, we observe 
	\begin{align*}
		L^D_{f_{\theta_1}}[\scrG^D(\theta_1)-\scrG^D(\theta_2)] =	L^D_{f_{\theta_1}}[u_{f_{\theta_1}}-u_{f_{\theta_2}}] = \divergence((f_{\theta_2}-f_{\theta_1})\cdot\nabla u_{f_{\theta_2}}),
	\end{align*}
	where the right-hand side defines a function in $C^\nu(\domain)$ for some small $\nu\in\pRZ$. \magenta{Thus}, applying the linear inverse $V_{f_{\theta_1}}^D$ on both sides of the equation, we obtain
	\begin{equation*}
	u_{f_{\theta_1}}-u_{f_{\theta_2}}= V_f^D[ \divergence((f_{\theta_2}-f_{\theta_1})\cdot\nabla u_{f_{\theta_2}})],
	\end{equation*}
	\magenta{and hence}, by an application of \cref{eq:darcy:inftynormbound}
	\begin{equation*}
		\norm{u_{f_{\theta_1}}-u_{f_{\theta_2}}}_\infty \leq \const(d,\domain,f_{\min})\cdot\norm{\divergence((f_{\theta_2}-f_{\theta_1})\cdot\nabla u_{f_{\theta_2}})}_\infty.
	\end{equation*}
	\magenta{Utilizing the product rule, the right-hand side can be further upper bounded by}
	\begin{align*}
		\norm{\divergence((f_{\theta_2}-f_{\theta_1})\cdot\nabla u_{f_{\theta_2}})}_\infty &\leq \norm{\nabla (f_{\theta_2}-f_{\theta_1})\cdot \nabla u_{f_{\theta_2}}}_\infty + \norm{(f_{\theta_2}-f_{\theta_1})\Delta u_{f_{\theta_2}}}_\infty\\
		%& \leq \norm{\nabla (f_{\theta_2}-f_{\theta_1})\cdot \nabla u_{f_{\theta_2}}}_\infty + \norm{f_{\theta_2}-f_{\theta_1}}_\infty\cdot \norm{\Delta u_{f_{\theta_2}}}_\infty\\
		& \leq \norm{f_{\theta_2}-f_{\theta_1}}_{C^1(\domain)}\cdot\left( \norm{u_{f_{\theta_2}}}_{C^1(\domain)}+ \norm{\Delta u_{f_{\theta_2}}}_\infty\right).
	\end{align*}
	We study the two remaining terms in the last display seperately, starting with the second one. We have 
		\begin{align*}
			\norm{\Delta u_{f_{\theta_2}}}_{\infty} & = \norm{f_{\theta_2}^{-1}\cdot(g-\nabla u_{f_{\theta_2}}\cdot \nabla f_{\theta_2})}_{\infty}\\
			&\leq \norm{f_{\theta_2}^{-1}}_{\infty}\cdot\left(\norm{g}_\infty + \norm{u_{f_{\theta_2}}}_{C^1(\domain)} \norm{ f_{\theta_2}}_{C^1(\domain)}\right)\\
			&\leq \const(\fmin)\cdot\left(\norm{g}_\infty + \norm{u_{f_{\theta_2}}}_{\Zygmund{2}{\domain}} \norm{ f_{\theta_2}}_{C^1(\domain)}\right),
		\end{align*}
		\magenta{where in the last step we have used that $f_{\theta_2}>\fmin$ and $\Zygmund{2}{\domain}\hookrightarrow C^1(\domain)$.}
		By an application of Lemma 22 in \cite{NicklGeerWang2020} and the continuous embedding $C^1(\domain)\hookrightarrow\Zygmund{1}{\domain}$, we \magenta{deduce}
		\begin{equation*}
			\norm{u_{f_{\theta_2}}}_{\Zygmund{2}{\domain}} \leq  \const(d,\domain,\norm{g}_\infty,f_{\min})\cdot (1+\norm{f_{\theta_2}}_{C^1(\domain)}^2)
		\end{equation*}
		\magenta{and hence}
		\begin{equation*}
			\norm{\Delta u_{f_{\theta_2}}}_{\infty} \leq \const(\fmin,\norm{g}_\infty,d,\domain)\cdot(1+\norm{f_{\theta_2}}_{C^1(\domain)}^3).
		\end{equation*} 
		\magenta{Another application of Lemma 22 in \cite{NicklGeerWang2020} yields}
		\begin{equation*}
			\norm{\divergence((f_{\theta_2}-f_{\theta_1})\cdot\nabla u_{f_{\theta_2}})}_\infty \leq \const(\fmin,\norm{g}_\infty,d,\domain)\cdot(1+\norm{f_{\theta_2}}_{C^1(\domain)}^3)\cdot  \norm{f_{\theta_2}-f_{\theta_1}}_{C^1(\domain)},
		\end{equation*}
		\magenta{which leads overall to}
		\begin{multline*}
			\norm{\scrG^D(\theta_1)-\scrG^D(\theta_2)}_\infty\\\leq \const(\fmin,\norm{g}_\infty,d,\domain)\cdot(1+\norm{f_{\theta_1}}_{C^1(\domain)}^3\lor\norm{f_{\theta_2}}_{C^1(\domain)}^3)  \norm{f_{\theta_2}-f_{\theta_1}}_{C^1(\domain)}.
		\end{multline*}
		\magenta{Now utilizing the properties of the regular link function $\Psi$ formulated in \cref{lem:linkfuncproperties} and Lemma 29 in \cite{NicklGeerWang2020}, it follows}
		\begin{multline*}
			\norm{\scrG^D(\theta_1)-\scrG^D(\theta_2)}_\infty\\
			 \leq  \const(\fmin,\norm{g}_\infty,d,\domain)\cdot(1+\norm{\theta_1}_{C^1(\domain)}^4\lor\norm{\theta_2}_{C^1(\domain)}^4)\cdot  \norm{\theta_2-\theta_1}_{C^1(\domain)}.
		\end{multline*}
		Hence, the forward map $\scrG^D$ satisfies C3) in \cref{ass:ForwardMapI} with $\gamma_2 = 4$, and $\eta_1=\eta_2 = 1$.
\end{enumerate}
After prooving the conditions C1) - C3) for the forward map $\scrG^D$ formulated in \cref{ass:ForwardMapI}, we are showing an analogous inverse continuity modulus for the solution operator $G^D$ of the Darcy Problem as assumed in C7) in \cref{ass:ForwardMapII}. The proof is based on a interpolation inequality for Sobolev spaces, which reads for all $\alpha\in\pRZ$, $\beta\in[0,\alpha+1]$
\begin{equation}\label{eq:Interpol1}
	\norm{u}_{H^\beta(\domain)} \leq \const(\alpha,\beta,\domain)\cdot\norm{u}_{\IL^2(\domain)}^{\frac{\alpha+1-\beta}{\alpha+1}}\cdot \norm{u}_{H^{\alpha+1}(\domain)}^{\frac{\beta}{\alpha+1}},
\end{equation}
\magenta{which is proven for instance in Theorem 1 in Section 4.3.1 in \cite{Triebel_1978} observing that $H^\beta(\domain)=B_{22}^\beta(\domain)$.}
\begin{enumerate}[label = {C'}7)]
	\item Inverse Continuity Modulus: Let $f_1,f_2\in\calF$ for some $\alpha>2+\frac{d}{2}$. \magenta{Assume there exists a constant $\operatorname{B}\in\pRZ$, such that $\norm{f_1}_{C^1(\domain)}\lor\norm{G^D(f_1)}_{C^2(\domain)}\leq\operatorname{B}$ and additionally that the heat source function $g$ satisfies $\inf_{x\in\domain}g(x) \geq g_{\min}$.} Lemma 24 in \cite{NicklGeerWang2020} \magenta{then} provides the existence of a finite constant $\const_{1}=\const_{1}(g_{\min},f_{\min},\operatorname{B},\domain,d)\in\pRZ$, such that 
	\begin{equation*}
		\norm{f_1-f_2}_{\IL^2(\domain)} \leq \const_{1}\cdot\norm{f_2}_{C^1(\domain)}\cdot\norm{G^{D}(f_1)-G^D(f_2)}_{H^2(\domain)}.
	\end{equation*}
	By an application of \cref{eq:Interpol1} \magenta{with $\beta =2$, we further have}
	\begin{multline*}
		\norm{G^{D}(f_1)-G^D(f_2)}_{H^2(\domain)}\\
		\leq \const_{2}(\alpha,\domain)\norm{G^{D}(f_1)-G^D(f_2)}_{H^{\alpha+1}(\domain)}^{\frac{2}{\alpha+1}}\cdot\norm{G^{D}(f_1)-G^D(f_2)}_{\IL^2(\domain)}^{\frac{\alpha-1}{\alpha+1}}.
	\end{multline*}
	Moreover, by an application of Lemma 23 in \cite{NicklGeerWang2020}, we have 
	\begin{equation*}
		\norm{G^{D}(f_1)-G^D(f_2)}_{H^{\alpha+1}(\domain)}^{\frac{2}{\alpha+1}} \leq \const_{3}\cdot(1+\norm{f_1}_{H^\alpha(\domain)}^{2\alpha}\lor\norm{f_2}_{H^\alpha(\domain)}^{2\alpha})
	\end{equation*}
	with a finite constant $\const_{3}=\const_{3}(\alpha,d,\domain,f_{\min},\norm{g}_{H^{\alpha-1}(\domain)})\in\pRZ$. Thus, we have 
	\begin{equation*}
			\norm{f_1-f_2}_{\IL^2(\domain)} \leq \const_{4}\cdot(1+\norm{f_1}_{H^\alpha(\domain)}^{2\alpha}\lor\norm{f_2}_{H^\alpha(\domain)}^{2\alpha})\cdot\norm{G^{D}(f_1)-G^D(f_2)}_{\IL^2(\domain)}^{\frac{\alpha-1}{\alpha+1}}
	\end{equation*}
	for a constant $\const_{4} = \const_{4}(\alpha,d,\domain,f_{\min},g_{\min},\operatorname{B},\norm{g}_{\alpha-1})\in\pRZ$.\\
	\magenta{Now, let $\operatorname{M}\in\pRZ$ be arbitrary and let $f_1$ and $f_2$ be chosen in such a way that $\norm{f_1}_{H^\alpha(\domain)}+\norm{f_2}_{H^\alpha(\domain)}\leq\operatorname{M}$ and $\norm{G^{D}(f_1)-G^D(f_2)}_{\IL^2(\domain)}\leq\delta$ for $\delta\in\pRZ$ small enough. By an application of Lemma 23 in \cite{NicklGeerWang2020} and the embedding $H^\alpha(\domain)\hookrightarrow C^1(\domain)$, we have $\norm{f_1}_{C^1(\domain)}\lor\norm{G^D(f_1)}_{C^2(\domain)}\leq c(\alpha,d,\domain,f_{\min},\norm{g}_{H^{\alpha-1}(\domain)})$. 
	By above computations we derive
	 \begin{equation*}
			\norm{f_1-f_2}_{\IL^2(\domain)} \leq \const_{5}\cdot\delta^{\frac{\alpha-1}{\alpha+1}}
	\end{equation*}
	with $\const_{5}=\const_{5}(\alpha,d,\domain,f_{\min},g_{\min},\operatorname{M},\norm{g}_{\alpha-1})\in\pRZ$. Taking the supremum over all admissible $f_1$ and $f_2$ as in the spirit of \cref{ass:ForwardMapII} shows the claim with $\tau = \frac{\alpha-1}{\alpha+1}\in(0,1)$.}
\end{enumerate}
We are now providing the proof of \cref{thm:main:darcy}.
\begin{proof}[Proof of \cref{thm:main:darcy}]	
	Firstly, note that the Darcy problem is formulated as special case of \cref{sec:model} with $V = W = \IR$. 
	\begin{enumerate}
		\item[i)] The existence of a maximizer $\fMAP$ follows directly by the equivalent formulation presented in \cref{sec:PDE} and an application of \cref{thm:maintheorem} under considerations that the forward map $\scrG^D$ satisfies all assumption of \cref{thm:maintheorem}.
		
		\item[ii)] As the forward map $\scrG^D$ satisfies C1) - C3) in \cref{ass:ForwardMapI}, we apply \cref{thm:maintheorem} directly with $\thetatilde = \theta_\frako$ by observing that
		\begin{equation*}
			\frakD_\rate^2(\fMAP,f_\frako) = \frakd_\rate^2(\thetaMAP,\theta_\frako).
		\end{equation*}
		Thus, for every $\bar{\const}\in\pRZ$ there exists a constant $c_1=c_1(\alpha,f_{\min},\norm{g}_\infty,d,\domain,\sigma)\in\pRZ$ sufficiently large, such that 
		\begin{equation*}
			\IP_{f_{\frako}}^N\left(\frakD_\rate^2(\fMAP,f_\frako)\geq\const_1(\rate_N^2\norm{\Psi^{-1}\circ f_\frako}_{H^\alpha(\domain)}^2+m^2)\right)\leq \const_2\exp\leq\left(-\bar{\const} m^2N\right)
		\end{equation*}
		for all $m\geq\rate_N$ satisfying \cref{eq:main:ass1} and \cref{eq:main:ass2}. We identify $f_\frako\in\calF$ with some $\theta_\frako = \Psi^{-1}\circ f_{\frako}\in H^{\alpha}_c(\domain)$ and conclude $\norm{\Psi^{-1}\circ f_\frako}_{H^\alpha(\domain)}\leq\const_{3}(\Psi,\operatorname{L})\in\pRZ$ uniformly over all $f_\frako\in\calF(\operatorname{L})$. Choose $\operatorname{M}\geq 1\lor\const_{3}$ and $m:=\operatorname{M}\rate_{N}\geq\rate_{N}$, such that $\const_1(\rate_N^2\norm{\Psi^{-1}\circ f_\frako}_{H^\alpha(\domain)}^2+m^2) \leq 2\const_1\operatorname{M}^2\rate_{N}^2$. Then, 
		\begin{equation*}
				\IP_{f_{\frako}}^N\left(\frakD_\rate^2(\fMAP,f_\frako)\geq2\const_1\operatorname{M}^2\rate_{N}^2\right)\leq \const_2\exp\leq\left(-\bar{\const} M^2 N\rate_N^2\right).
		\end{equation*}
		Updating constants yields the claim.
		\item[iii)] As we have shown an analogoue inverse continuity modulus C'7), the claim follows directly following along the lines of the proof of \cref{cor:stability} and is hence omitted. 
	\end{enumerate}
	\mbox{}
\end{proof}
\phantom{.}\\
We are now providing the proof of \cref{thm:main2:darcy}.
\begin{proof}[Proof of \cref{thm:main2:darcy}]
		\magenta{The proof follows the lines of the proofs of analogous results in the White Noise Model setting derived in \cite{NicklGeerWang2020} (see Proof of Theorem 9). Hence, in this work we only provide a sketch of the proof and ommit the details, which are left for the interested reader.}
		\begin{enumerate}
			\item[i)] By an application of \cref{eq:Interpol1}, we obtain 
			\begin{multline*}
				\norm{G^D(\fMAP)-G^D(f_\frako)}_{H^\beta(\domain)}^2\\
				 \leq \const_{1}\cdot \norm{G^D(\fMAP)-G^D(f_\frako)}_{\IL^2(\domain)}^{\frac{2(\alpha+1-\beta)}{\alpha+1}}\cdot\norm{G^D(\fMAP)-G^D(f_\frako)}_{H^{\alpha+1}(\domain)}^{\frac{2\beta}{\alpha+1}}.
			\end{multline*}
			\magenta{The first factor can be directly upper bounded by $\frakD_\rate^2(\thetaMAP,\theta_\frako)$. The second factor can be further upper bounded by firstly applying Lemma 23 in \cite{NicklGeerWang2020} and then by considering the definition of $\frakD_\rate^2(\thetaMAP,\theta_\frako)$.}
			Decomposing the underlying probability space into slices 
			\begin{equation*}
				A_0 :=\left\{\frakD_\rate^2(\fMAP,f_\frako)<\delta_N^2\right\}\text{ and } A_j:=\left\{\frakD_\rate^2(\fMAP,f_\frako)\in(2^{2(j-1)}\delta_N^2,2^{2j}\delta_N^2]\right\}
			\end{equation*}
			for some $\delta_N \propto \rate_{N}$, the corresponding expectation can be upper bounded by 
			\magenta{
			\begin{multline*}
				\IE_{f_\frako}^N\left[\norm{G^{D}(\fMAP)-G^D(f_\frako)}_{H^\beta(\domain)}^2\right]\\
				  \lesssim \sum_{j\in\IN_0}\IE_{f_\frako}^N\left[\mathds{1}_{A_j}\frakD_\rate^2(\thetaMAP,\theta_\frako)^{\frac{\alpha+1-\beta}{\alpha+1}}\left(1+(\rate_{N}^{-2}\frakD_\rate^2(\thetaMAP,\theta_\frako))^{\alpha^2\beta}\right)\right]\\
				    \lesssim\delta_N^{\frac{2(\alpha+1-\beta)}{\alpha+1}}\bigg(1+ \sum_{j\in\IN_0}2^{2j\frac{\alpha+1-\beta}{\alpha+1}}(1+2^{2j\alpha^2\beta})\IP_{f_\frako}^N(A_j)\bigg).
			\end{multline*}}
			Using the concentration inequality of \cref{thm:main:darcy}, we obtain
			\begin{equation*}
				\IE_{f_\frako}^N\left[\norm{G^{D}(\fMAP)-G^D(f_\frako)}_{H^\beta(\domain)}^2\right] \lesssim \delta_N^{\frac{2(\alpha+1-\beta)}{\alpha+1}}\cdot(1+o(N^{-\frac{1}{2}})),
			\end{equation*}
			which shows the claim.
			
			\item[ii)] Similiarly, as in the proof of C'7) before, we upper bound 
			\begin{equation*}
				\norm{\fMAP-f_\frako}_{\IL^2(\domain)}^2 \lesssim\norm{\fMAP}_{C^1(\domain)}^2\cdot\norm{G^D(\fMAP)-G^D(f_\frako)}_{H^2(\domain)}^2.
			\end{equation*}
			Thus, by an application of \cref{eq:Interpol1} with $\beta=2$, we obtain
				\begin{multline*}
				\norm{G^D(\fMAP-G^D(f_\frako)}_{H^{2}(\domain)}^2\\ \lesssim\norm{G^D(\fMAP-G^D(f_\frako)}_{H^{1+\alpha}(\domain)}^{\frac{4}{1+\alpha}} 	\norm{G^D(\fMAP)-G^D(f_\frako)}_{\IL^{2}(\domain)}^{\frac{2(\alpha-1)}{\alpha+1}}.
			\end{multline*}
			\magenta{Using the same techniques as for the prediction error, we obtain the bound}
			\begin{equation*}
				\IE_{f_\frako}^N\left[\norm{\fMAP-f_\frako}_{\IL^2(\domain)}^2\right] \lesssim\delta_N^{\frac{2(\alpha-1)}{\alpha+1}}\cdot(1+o(N^{-\frac{1}{2}})),
			\end{equation*}
			which proves the claim.
		\end{enumerate}
		\mbox{}
\end{proof}

%% file: content/appendix/auxres.tex
\subsection{Auxiliary results}

The following lemma was proven in \cite{NicklWang2023} (see Lemma 3.12) and is a key tool in the proof of \cref{thm:maintheorem}.
\begin{lemma}[Chaining for an Empirical Process]\label{lem:EPT}
	Let $\Theta$ be a countable set. Let $(\Omega,\calA,\IP^Z)$ be a probability space. Suppose a class of real-valued measurable functions
	\begin{equation*}
		\calH :=\left\{h_\theta:\Omega\to\IR:\,\theta\in\Theta\right\}.
	\end{equation*}
	Assume that $\calH$ is uniformly bounded and has a variance envelope, namely
	\begin{enumerate}
		\item[i)] there exists $\Const_{\operatorname{B}}\in\pRZ$, such that 
			\begin{equation*}
				\sup_{\theta\in\Theta}\left\{\sup_{\omega\in\Omega}\left\{|h_\theta(\omega)|\right\}\right\} \leq  \Const_{\operatorname{B}},
			\end{equation*}
		\item[ii)] there exists  $\sigma^2\in\pRZ$, such that
		\begin{equation*}
			\sup_{\theta\in\Theta}\left\{\IE^Z [h_\theta(Z)^2]\right\}\leq \sigma^2,
		\end{equation*}
		 for some $Z\sim\IP^Z$, where $\IE^Z$ denotes the expectation operator associated to $\IP^Z$.
	\end{enumerate}
	We define the metric entropy integrals 
	\begin{equation*}
		\calJ_2(\calH) := \int_0^{4\sigma}\sqrt{\log N(\calH,d_2,\rho)}\mathrm{d}\rho\text{ and }	\calJ_\infty(\calH) := \int_0^{4\Const_{\operatorname{B}}}\log N(\calH,d_\infty,\rho)\mathrm{d}\rho
	\end{equation*}
	with $\IL^2$- and $\IL^\infty$-structure given by 
	\begin{equation*}
		d_2(\theta,\theta'):= \sqrt{\IE^Z[h_\theta(Z)-h_{\theta'}(Z)]^2}\text{ and }d_\infty(\theta,\theta'):=\sup_{\omega\in\Omega}\left\{|h_\theta(\omega)-h_{\theta'}(\omega)|\right\},
	\end{equation*}
	respectively. Let $N\in\IN$. For i.i.d. copies $Z_1,\dots,Z_N$ of $Z\sim \IP^Z$ and $\varepsilon_1,\dots\varepsilon_N$ of $\varepsilon\sim\operatorname{N}(0,1)$ independent of $Z$, consider empirical processes either of the form
	\begin{equation*}
		T_N(\theta):=\frac{1}{\sqrt{N}}\sum_{i\in\nset{N}}h_\theta(Z_i)\varepsilon_i, \, \theta\in\Theta
	\end{equation*}
	or 
	\begin{equation*}
		T_N(\theta):=\frac{1}{\sqrt{N}}\sum_{i\in\nset{N}}(h_\theta(Z_i)-\IE[h_\theta(Z)]),\, \theta\in\Theta.
	\end{equation*}
	Then, there exists a universal constant $L\in\pRZ$, such that for all $x\in\IR_{\geq 1}$, we have
	\begin{equation*}
		\IP\left(\sup_{\theta\in\Theta}\left\{|T_N(\theta)|\right\}\geq L\left(\calJ_2(\calH)+\sigma\sqrt{x}+N^{-\frac{1}{2}}(\calJ_\infty(\calH)+\Const_{\operatorname{B}}x)\right)\right) \leq 2\exp(-x).
	\end{equation*}
\end{lemma}
Similarly to \cite{NicklWang2023}, whenever we apply \cref{lem:EPT} in this work with an uncountable set $\Theta$, one can verify that the supremum can be realized as one over a countable subset of it.

%In the proof of \cref{thm:maintheorem}, we need the following lemma providing an upper bound for covering numbers in finite product spaces.

\begin{lemma}[Covering Numbers of Product Spaces]\label{lem:CoveringPS}
	Let $(\tilde{S},\norm{\cdot}_{\tilde{S}})$ be a normed space. Assume there exists a function $\psi:\pRZ\times\pRZ\to\IR$, such that for any $r,\rho\in\pRZ$ the covering number $N(\IB_{\tilde{S}}(0,r),\norm{\cdot}_{\tilde{S}},\rho) = \psi(r,\rho)$ is finite. For $K\in\IN$ define $S:=\bigtimes_{i\in\nset{K}} \tilde{S}$ equipped with the norm $\norm{s}_S := \left(\sum_{i\in\nset{K}}\norm{s_i}_{\tilde{S}}^2\right)^{\frac{1}{2}}$	for any $s = (s_1,\dots,s_K)\in S$. We then have 
	\begin{equation*}
		N(\IB_{S}(0,r),\norm{\cdot}_S,\rho)  \leq \psi\left(r,\frac{\rho}{\sqrt{K}}\right)^K.
	\end{equation*}
\end{lemma}
\begin{proof}[Proof of \cref{lem:CoveringPS}]
	By hypothesis, for all $\rho,r\in\pRZ$ there exists a covering $\calC_{\tilde{S}}:=\{\tilde{s}^{(i)}\}_{i\in\nset{N}}\subseteq \tilde{S}$ of $\IB_{\tilde{S}}(0,r)$ with $N:=\psi(r,\rho)$, i.e. for all $\tilde{s}\in\IB_{\tilde{S}}(0,r)$ there exists $i_\circ=i_\circ(\tilde{s})\in\nset{N}$ with $\norm{\tilde{s}-\tilde{s}^{(i_\circ)}}_{\tilde{S}} \leq \rho$.
	Now let $s=(s_1,\dots,s_K)\in\IB_S(0,r)$ be arbitrary chosen, such that
	\begin{equation*}
		\norm{s}_S^2 = \sum_{j\in\nset{K}}\norm{s_j}_{\tilde{S}}^2 \leq r^2,
	\end{equation*}
	and thus, we have $s_j\in \IB_{\tilde{S}}(0,r)$ for all $j\in\nset{K}$. Hence, for each $s_j$, $j\in\nset{K}$, there exists $i=i_j\in\nset{N}$, such that $s_j\in\IB_{\tilde{S}}(\tilde{s}^{(i_j)},\rho)$.
	Defining $s^{(i_1,\dots,i_K)}:= (\tilde{s}^{(i_1)},\dots,\tilde{s}^{(i_K)})\in S$, we have 
	\begin{equation*}
		\norm{s-s^{(i_1,\dots,i_K)}}_S^2 = \sum_{j\in\nset{K}}\norm{s_j-\tilde{s}^{(i_j)}}_{\tilde{S}}^2 < K\rho^2.
	\end{equation*}
	Hence, $\calC_S := \bigtimes_{i\in\nset{K}}\calC_{\tilde{S}}$ forms a $\sqrt{K}\varepsilon$-covering of $\IB_S(0,r)$ and thus
	\begin{equation*}
		N(\IB_{S}(0,r),\norm{\cdot}_S,\sqrt{K}\rho) \leq N^K = \psi(r,\rho)^K.
	\end{equation*}
	This shows the claim.
\end{proof}

In the following lemma we use \cref{lem:CoveringPS} to derive an upper bound for the covering number of Sobolev balls.

\begin{lemma}[Covering Numbers of Sobolev Spaces]\label{lem:SobolevCovering}
	Let $d\in\IN$. Let $\domain\subseteq\IR^d$ be a bounded domain with smooth boundary $\partial\domain$. Let $(W,\Wnorm{\cdot})$ be a finite-dimensional normed vector space with $\operatorname{dim}_\IR W = d_W$. 
	\begin{enumerate}
		\item[i)] Let $\kappa\in\pRZz$, $r\in\pRZ$.
		Assume that $\calF$ is a subspace of $H^\alpha(\domain,W,r)$, if $\kappa<\frac{1}{2}$ and of $H^\alpha_c(\domain,W,r)$, otherwise. Then there exists a constant $c=c(\alpha,\kappa,d)\in\pRZ$, such that 
		\begin{equation*}
			N(\calF,\norm{\cdot}_{(H^\kappa(\domain,W))^*},\rho) \leq \exp\left(\const d_W \left( r\sqrt{d_W}\rho^{-1}\right)^{\frac{d}{\alpha+\kappa}}\right)
		\end{equation*}
		for any $\rho\in\pRZ$.
		\item[ii)] Let $\eta\in\pRZz$ , $r\in\pRZ$. Assume that $\calF$ is a subspace of $H^\alpha(\domain,W,r)$ for $\alpha > \eta + \frac{d}{2}$. Then there exists a constant $c=c(\alpha,\eta,d)\in\pRZ$, such that 
		\begin{equation*}
			N(\calF,\norm{\cdot}_{C^{\eta}(\domain,W)},\rho) \leq \exp\left(\const d_W \left( r\sqrt{d_W}\rho^{-1}\right)^{\frac{d}{\alpha-\eta}}\right)
		\end{equation*}
		for any $\rho\in\pRZ$.
	\end{enumerate}
\end{lemma}
\begin{proof}[Proof of \cref{lem:SobolevCovering}]
	\mbox{}
	\begin{enumerate}[label = i)]
		\item[i)] Restricting ourselfes to the case $d_W=1$, the claim is already shown \magenta{ in the proof of Lemma 19} in \cite{NicklGeerWang2020}, \magenta{where a differentiation between the cases $\kappa\geq\frac{1}{2}$ and $\kappa<\frac{1}{2}$ is necessary to find proper $(H^\kappa(\domain))^*$-coverings of $H^\alpha_c(\domain)$ and $H^\alpha(\domain)$, respectively.} For general $d_W\in\IN$, the claim follows by an application of \cref{lem:CoveringPS}.
		\item[ii)]
			We restrict ourselves to the case $d_W = 1$. The general claim follows then by \cref{lem:CoveringPS}. First, let $\eta = 0$, such that $\norm{\cdot}_{C^\eta(\domain)}=\norm{\cdot}_\infty$. \magenta{For $d=1$, if $\domain$ equals $[0,1]$, in \cite{Gine2016} a covering result is shown, see Corollary 4.3.38 for instance.}
			Using usual extension arguments, one obtains \magenta{for the general case}
			\begin{equation*}
				N(\calF,\norm{\cdot}_{\infty},\rho) \leq \exp\left(\const(\alpha,d)\cdot\left(\frac{r}{\rho}\right)^{\frac{d}{\alpha}}\right),
			\end{equation*}
			which is also used in \cite{Nickl_2023}. Now, let $\eta\in\pRZ\setminus\IN$. Following \cref{rem:weakconvergence}, the embedding $H^\alpha(\domain)\hookrightarrow C^{\eta}(\domain) = B_{\infty\infty}^{\eta}(\domain)$ is compact, and by \cite{Triebel1996} \magenta{(Section 3.3) and \cite{Triebel_1978} (Section 4.10.3), respectively}, we further have
			\begin{equation*}
				N(\calF,\norm{\cdot}_{C^{\eta}(\domain)},\rho) \leq \exp\left(\const\cdot\left(\frac{r}{\rho}\right)^{\frac{d}{\alpha-\eta}}\right).
			\end{equation*}
			Now let $\eta = k\in\IN$. According to \cref{rem:weakconvergence}, the embedding $H^\alpha(\domain)\hookrightarrow B_{\infty1}^k(\domain)\hookrightarrow C^k(\domain)$ is compact, as $\alpha > k +\frac{d}{2}$. Thus, for any $\theta_1,\theta_2\in H^\alpha(\domain)$, we have
			\begin{equation*}
				\norm{\theta_1-\theta_2}_{C^{\eta}(\domain)} \leq \const_{\texttt{em}} \cdot \norm{\theta_1-\theta_2}_{B_{\infty1}^k(\domain)}.
			\end{equation*}
			Thus, it sufficies to find an covering of $H^\alpha(\domain,r)$ w.r.t the $\norm{\cdot}_{B_{\infty1}^k(\domain)}$-norm. Following the theory in \cite{Triebel1996} \magenta{(Section 3.3) and \cite{Triebel_1978} (Section 4.10.3), respectively}, we have 
			\begin{equation*}
				N(H^\alpha(\domain,r),\norm{\cdot}_{B_{\infty1}^k(\domain)},\rho) \leq \exp\left(\const(\alpha,k,d)\cdot\left(\frac{r}{\rho}\right)^{\frac{d}{\alpha-k}}\right),
			\end{equation*}
			and hence
			\begin{equation*}
				N(H^\alpha(\domain,r),\norm{\cdot}_{C^k(\domain)},\rho) \leq \exp\left(\const(\alpha,k,d)\cdot\left(\frac{\const_{\texttt{em}}r}{\rho}\right)^{\frac{d}{\alpha-k}}\right).
			\end{equation*}
			The claim follows.
	\end{enumerate}
\end{proof}

%In the following lemma, we summarize all properties of regular link function as defined in \cref{def:LinkFunction}, which we frequently use in the present work.
\begin{lemma}[Properties of Regular Link Functions]\label{lem:linkfuncproperties}
	Let $\Psi$ be a regular link function as in \cref{def:LinkFunction}. We have the following properties.
	\begin{enumerate}
		\item[i)] Let $\theta_1,\theta_2\in C^0(\domain)$. Then \magenta{ there exists a constant $c=c(\Psi)\in\IR_{>0}$, such that}
		\begin{equation*}
			\norm{\Psi\circ\theta_1-\Psi\circ\theta_2}_{\infty} \leq \const(\Psi)\cdot \norm{\theta_1-\theta_2}_\infty.
		\end{equation*}
			\item[ii)] Let $\theta_1,\theta_2\in C^1(\domain)$. Then \magenta{ there exists a constant $c=c(\Psi)\in\IR_{>0}$, such that}
		\begin{equation*}
			\norm{\Psi\circ\theta_1-\Psi\circ\theta_2}_{C^1(\domain)} \leq \const(\Psi)\cdot(1+\norm{\theta_1}_{C^1(\domain)}\lor \norm{\theta_2}_{C^1(\domain)})\cdot \norm{\theta_1-\theta_2}_{C^1(\domain)}.
		\end{equation*}
			\item[iii)] Let $\theta_1,\theta_2\in \IL^2(\domain)$. Then \magenta{ there exists a constant $c=c(\Psi)\in\IR_{>0}$, such that}
		\begin{equation*}
			\norm{\Psi\circ\theta_1-\Psi\circ\theta_2}_{\IL^2(\domain)} \leq \const(\Psi)\cdot \norm{\theta_1-\theta_2}_{\IL^2(\domain)}.
		\end{equation*}
	\end{enumerate}
\end{lemma}
\begin{proof}[Proof of \cref{lem:linkfuncproperties}]
	Elementary computations show all three statements of the Lemma, since by hypothesis $\Psi$ as well as all its derivatives are bounded and hence Lipschitz continuous.
\end{proof}

\let\cleardoublepage\clearpage